\documentclass[11pt]{article} 
\usepackage{amsmath, amssymb} 
\usepackage{graphicx}
 \usepackage{graphics}

\newcommand{\rit}{\mathbb R}

\newcommand{\cit}{\mathbb C}      
\newcommand{\zit}{\mathbb Z}
\newcommand{\nin}{\notin}

\newcommand{\eps}{\varepsilon}

\newcommand{\alali}{ \mbox{ } \\}
\newcommand{\dis}{\displaystyle}
\newcommand{\ra}{\rightarrow}

\newcommand{\bqr}{\begin{eqnarray}}
\newcommand{\bqre}{\begin{eqnarray*}}
\newcommand{\eqr}{\end{eqnarray}}
\newcommand{\eqre}{\end{eqnarray*}}

\newcommand{\cqfd}{{\nobreak\hfil\penalty50\hskip2em\hbox{}\nobreak\hfil
$\square$\qquad\parfillskip=0pt\finalhyphendemerits=0\par\medskip}}

\newcommand{\bpro}{\noindent{\rm\bf Proof : }}
\newcommand{\epro} {\cqfd} 

\newtheorem{theorem}{Theorem}
\newtheorem{lemma}{Lemma}
\newtheorem{proposition}{Proposition}
\newtheorem{corollary}{Corollary}

\newtheorem{remark}{Remark}

\renewcommand{\thefigure}
                             {\arabic{figure}}
 \renewcommand{\theequation}
                             {\arabic{equation}}

\title{   The Method of Strained Coordinates 
          for Vibrations
         with Weak  Unilateral Springs  }

  \author{   St\'ephane Junca 
             \thanks{ Universit\'{e} de Nice Sophia-Antipolis,  JAD laboratory, 
                               Parc Valrose, 06108 Nice, France,  junca@unice.fr}
          \;  \&  \; 
           Bernard Rousselet 
                \thanks{ Universit\'{e} de Nice Sophia-Antipolis,  JAD laboratory,
            Parc Valrose, 06108 Nice, France, br@unice.fr}
          }
\date{\empty}
   
\begin{document}

 \maketitle

  \abstract{We study some spring mass  models for a structure 
   having  some  unilateral springs  of small rigidity   
    $\eps$.
 We obtain and justify mathematically 
  an asymptotic expansion with the method of strained coordinates
 with new tools to handle such defects,
 including a non negligible cumulative effect over a long time: 
$T_\eps  \sim 1/\eps$ as usual;  or, for a new critical case,
 we can only expect:  $T_\eps \sim 1/\sqrt{\eps}$.
We check numerically these results and present a purely numerical algorithm to compute ``Non linear Normal Modes'' (NNM); this algorithm provides results close to the asymptotic expansions but enables us to compute NNM even when $\eps$ becomes larger. 
  }  \bigskip \\
   {\bf Keywords}: 
   nonlinear vibrations, 
  method of 
 strained coordinates,  
 piecewise linear, 
 unilateral spring, 
 approximate nonlinear normal mode.
             \\

 $\mbox{ }$\hspace{-6mm} 
{\it Mathematics Subject Classification.}
  Primary:  34E15; \\
$\mbox{ }$ \hspace{53mm} Secondary: 26A16, 26A45, 41A80. 
 
       

    \section{Introduction \label{SI}} 

 For spring mass  models, the presence of 
 a small piecewise linear rigidity can model a small defect
 which implies unilateral reactions on the structure.
 So, the  nonlinear and 
 piecewise linear function $u_+=\max(0,u)$ plays a key role in this paper.
   For nondestructive testing we  study a  non-smooth nonlinear effect
 for large time by asymptotic expansion of the vibrations.
 New features and comparisons with classical  cases 
of smooth perturbations are given, for instance, with  the classical
 Duffing equation:
 $ \ddot{u}+u +\eps u^3=0$ and the non classical case:
  $ \ddot{u}+u +\eps u_+=0$.
Indeed, piecewise  linearity is non-smooth:  nonlinear 
 and Lipschitz but not differentiable. 
 We give some new results to validate such asymptotic expansions.
 Furthermore, these tools are also valid 
for a more general   non linearity.
 A nonlinear crack approach for elastic waves 
 can be found in \cite{JuLoA}.  
  Another approach   in the framework of  non-smooth analysis  
 can be found in \cite{Att,CHL,PaoSch}.
\\
    For short time, a linearization procedure is enough 
  to compute a good approximation.  
   But for large time, nonlinear cumulative effects 
  drastically alter the nature of the solution. 
    We will consider the classical 
  method of strained coordinates  to compute asymptotic expansions.
 The idea goes further back to Stokes, who in 1847 calculated 
 periodic solutions for a weakly nonlinear wave propagation problem,
 see \cite{KC68,KC96,Mi,Nayfeh} for more details and references therein. 
 Subsequent authors have generally referred to this as the 
  Poincar\'e method  or the Lindstedt method.
 It is a simple and efficient method which gives us
  approximate nonlinear normal modes
 with $1$ or more degrees of freedom.
\\
 Lindstedt-Poincar\'e method has been already used in \cite{Ves}
 to study  NNM of a piecewise linear system with two degrees of freedom. 
 Here the non linearity is somewhat more general. We consider $N$ dimensional systems. Moreover we prove rigorously  the validity of the expansion. On  the other hand \cite{Ves} addresses other very
interesting open problems such as: bifurcation of solutions, higher order expansions, stability of solutions.
 \\
   In section \ref{s1dof} we present the method on 
  an explicit case with  an internal Lipschitz force. 
We  focus on an  equation with   one degree of freedom with expansions 
 valid for time   of order $\eps^{-1}$
 or, more surprisingly, $ \eps^{-1/2}$ for a degenerate contact.
  \\
   Section \ref{sexp+} contains a tool to expand $(u+\eps v)_+$
 and some accurate estimate for the remainder. 
 This is a new key point  to validate
 the  method of strained coordinates with unilateral contact.
 \\  
  In Section \ref{sNdof}, 
  we extend previous results to  systems 
 with $N$ degrees of freedom,
first, 
with the same accuracy  for approximate nonlinear normal modes,
 then, with less accuracy with all modes.
We check numerically these results and present a purely numerical algorithm to compute ``Non linear normal Modes'' (NNM) in the sense of Rosenberg \cite{Rosenberg}; see \cite{ABBC06} for two methods for the computation of NNM; see \cite{JPS} for a computation of non linear normal mode with unilateral contact and \cite{nnm-kpgv} for a synthesis on non linear normal modes;
 this algorithm provides results close to the asymptotic expansions but enables to compute NNM even when $\eps$ becomes larger. 
\\
In Section \ref{scos}, we briefly  explain why we only perform expansions with 
 even periodic functions to compute the nonlinear frequency shift.
\\
   Section \ref{sA} is an appendix 
 containing some technical proofs and results.

    \section{One degree of freedom \label{s1dof}}  

    \subsection{Explicit  angular frequency   \label{ss1dofexp}}  
\begin{figure}[!ht]
\begin{center}
\includegraphics[angle=0, width=11cm, height=3cm]{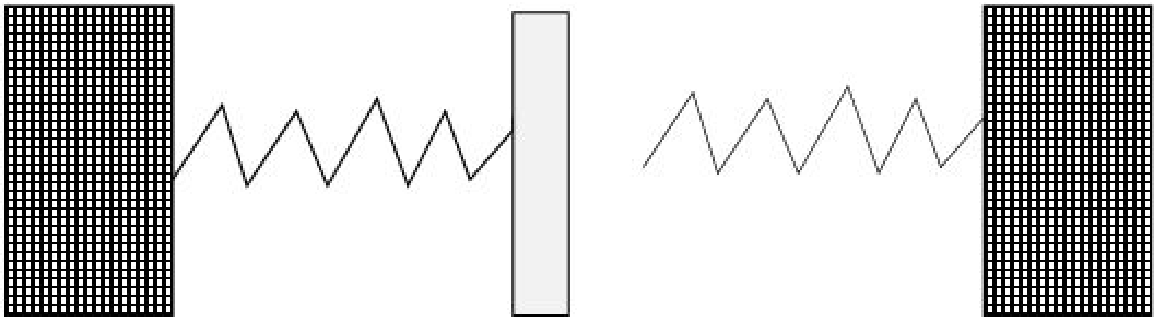}
\caption{ Two springs,   on the right
    it has only  a unilateral contact. \label{Ressort}}
\end{center}
\end{figure}

We consider a one degree of freedom spring-mass system (see figure \ref{Ressort}): one spring is  classical linear and attached to the mass and to a rigid wall, the second is still linear attached to a rigid wall  but has a unilateral contact with the mass; this is to be considered as a damaged spring. 
The force acting on the mass is  $k_1 u + k_2 u_+ $ where $u$ is the displacement of the mass $m$, $k_1$, the rigidity of the undamaged spring and $k_2$, the rigidity of the damaged unilateral spring. We notice that the 
 term 
 $ u_+$
is Lipschitz but not differentiable with respect to $u$. 
Assuming that $k_2=\widetilde{\eps} k_1$,  
 $\eps = \widetilde{\eps} \omega_0^2$  with $\omega_0^2=k_1/m$,
  we can consider the equation:
\bqr \label{ode}
 \ddot{u} + \omega_0^2 u + \eps u_+ =0,  
  \qquad \mbox{ with } u_+ = \max(0,u).
\eqr
The  associated  energy 
 is
 $  E= (\dot{u}^2 +\omega_0^2 u^2 + \eps(u_+)^2)/2.$
 Therefore, the level sets of $E(u,\dot{u})$ will
 be made of two half ellipses. 
 Indeed, for $u<0$ the level set is an half ellipse,
  and for $u>0$ is another  half ellipse.
  Any solution $u(t)$  is confined to a closed level curve 
  of $E(u,\dot{u})$ and is necessarily a periodic function of $t$.
 \\
  More precisely, a non trivial solution ($E>0$)
 is on the half  ellipse: $\dot{u}^2 + \omega_0^2 u =2 E$,
 in the phase plane
 during the time $T_C= \pi/\omega_0$,
 and on the half ellipse 
   $\dot{u}^2 + (\omega_0^2+\eps) u = 2E$
during the time  
 $T_E= \pi/\sqrt{\omega_0^2+\eps}$.
 Then the  period  is exactly 
 $ P(\eps) = \dis (1 + \left(1+\eps/\omega_0^2\right)^{-1/2})
            \pi/\omega_0,$
 and the exact angular frequency  is: 
\bqr \label{exactomega1}
 \dis 
  \omega(\eps)& =& 
        \dis  2 \omega_0 \left(1+ \left(1+\eps/\omega_0^2\right)^{-1/2}\right)^{-1} 
            = 
   \dis  \omega_0 + \frac{\eps}{(4\omega_0)} -\frac{\eps^2}{(8\omega_0^3)} 
                 + \mathcal{O}(\eps^3).
\eqr
Let us compare  with 
the angular frequency  $\omega_D(\eps)$ for Duffing
 equation
 where  the nonlinear term is  $u^3$ instead of $u_+$.
 $\omega_D(\eps)$ depends on the amplitude $a_0$
 of the solution (see for instance \cite{KC68,KC96,Mi,Nayfeh}):
$\dis
  \omega_D(\eps)   =  \omega_0 + \frac{ 3}{8\omega_0^2} a_0^2\eps
                    -\frac{ 15}{256\omega_0^4} a_0^4\eps^2
                   + \mathcal{O}(\eps^3).
$
%
    \subsection{The method of strained coordinates
             \label{ss1dofmsc}}  
\begin{figure}[!ht]
\begin{center}
\includegraphics[angle=0, width=12cm, height=7cm]{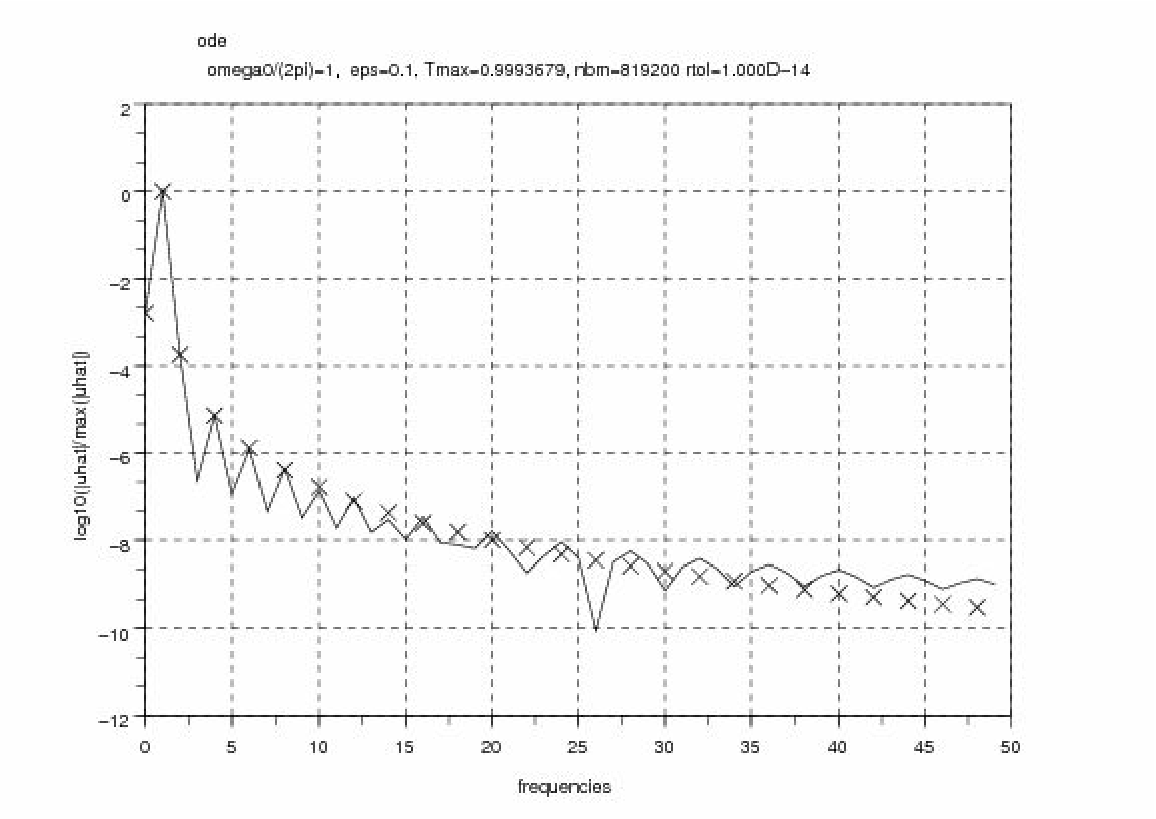}
\caption{ \label{figspectra}
 $v_1$,  
 $ \dis \log_{10} \left(\frac{|\widehat{u}_\eps|}{\max|\widehat{u}_\eps|} \right)  $
}
\end{center}
\end{figure}

%
 Now, we compute, with  the method of strained coordinates,
  $\omega_\eps$, {\it an approximation}
  of the exact angular frequency  $\omega(\eps)$ which is smooth 
 with respect to $\eps$ by exact formula
  (\ref{exactomega1}):
  $ 
       \omega(\eps)=\omega_\eps~+~\mbox{\rm O}(\eps^3).
 $ 
 We expound  this case completely to use the  same method
 of strained coordinates  later 
    when we will not have such  an explicit formula.
\\
Let us define the new time $ s  = \omega_\eps t $ 
and rewrite equation  (\ref{ode})
 with  $ v_\eps(s)   =  u_\eps(t) $  
  \bqr 
  \label{odes} 
\omega_\eps^2 v_\eps''(s) + \omega_0^2 v_\eps(s) + \eps (v_\eps(s))_+ =0,
&  \mbox{ with } s  = \omega_\eps t, 
       \qquad  u_\eps(t)   =  v_\eps(s),
 \eqr
To simplify, 
   $u_\eps$ is  subjected to  the following initial conditions: 
$
  u_\eps(0) = a_0 >0,\; \dot{u}_\eps(0) =0,$ 
i.e.  $v_\eps(0)=a_0$ and $v_\eps'(0)=0$.
 Similar computations are valid for negative $a_0$, see Proposition 
 \ref{propexplicit} below.  With more general data, 
i.e. when $\dot{u}_\eps(0) \neq 0$, computations are more complicate
 and give the same approximate angular frequency $\omega_\eps$, 
 see section \ref{scos}. 
\\
 In the  new time $s$, we use  the following {\bf ansatz} 
 \bqr  \label{ansatz} 
   \omega_\eps  = \omega_0+\eps \omega_1+\eps^2 \omega_2,
  \qquad  v_\eps(s)
        =  v_0(s)+\eps v_1(s ) +\eps^2 r_\eps(s).
 \eqr 
   $\omega_1$ and $ \omega_2$  are unknown. Since  $   \omega_\eps^2  =  
        \alpha_0+\eps \alpha_1+\eps^2 \alpha_2 +\mathcal{O}(\eps^3), \;
  %
  \alpha_0  = \omega_0^2,\; \alpha_1 = 2\omega_0\omega_1 ,
  \quad \alpha_2 = \omega_1^2+2\omega_0\omega_2,
  $
 we have to find   $\alpha_1$ and $ \alpha_2$.
\\
We will   also use the  following expansion,  
$ (u+\eps v)_+ =  u_+ + \eps H(u)v +  \eps\chi_\eps(u,v),
$ 
justified later, 
where $ H(.)$ is the Heaviside function, equal to $1$ if $u>0$ and else $0$.
Since $H(.)$ is not differentiable at  $u=0$,
 the remainder $\eps\chi_\eps(u,v)$
is not the classical Taylor's remainder. 
This lack of smoothness is a problem to validate mathematically  
the Lindstedt-Poincar\'e method.  
The remainder problem  is studied in section \ref{sexp+}  below.
\\
 Now, replacing   ansatz \eqref{ansatz} in  (\ref{odes}) 
   we obtain   differential
equations and initial data   for $v_0,v_1,r_\eps$ 
with $L(v) = -\alpha_0( v''  + v) $:
\bqr  
\label{eq0onedof}
 L(v_0) = & -\alpha_0( v''  + v) =0  , 
            &  v_0(0)=a_0,  v'_0(0)=0,\\
  \label{eq1onedof}
L( v_1)  = &  (v_0)_+ +  \alpha_1 v_0'', 
               & v_1(0)=0, \; \; v'_1(0)=0, \\       
 \label{rode}
L( r_\eps)  = & 
  H(v_0)v_1 + \alpha_2 v_0'' + \alpha_1 v_1'' 
 +  R_\eps(s), 
      &      r_\eps(0)=0,\; \; r'_\eps(0)=0.
\eqr 
  Now we  compute, 
  $\alpha_1$, $v_1$ and then $\alpha_2$.
  We have $v_0(s)= a_0 \cos(s)$. A key point in the method of
  strained coordinates is to keep bounded $v_1$ and $r_\eps$ for large time
 by a  choice of $\alpha_1$ for $u_1$ and $\alpha_2$ for $r_\eps$.
 For this purpose, we avoid resonant or {\bf secular} term 
in the right-hand-side of 
 equations (\ref{eq1onedof}), (\ref{rode}).
Let us first focus on $\alpha_1$. Notice that,
 $v_0(s)= a_0 \cos(s)$ and $a_0 > 0$, so 
 $ 
      \dis (v_0)_+   =  \dis  
   a_0 \left( \frac{\cos s}{2}+  \frac{|\cos s|}{2}\right).
$ 
  Note that $|\cos(s)|$ has  no term with frequencies $\pm 1$, 
     since there are  only even frequencies.
 Thus $( (v_0)_+ - \alpha_1 v_0)= 
  a_0\cos(s)(1/2 -\alpha_1) + a_0|\cos(s)|/2$ 
 has no secular term if and only if $\alpha_1 = 1/2$,  
 so $ \omega_1= 1/(4\omega_0)$. 
Now, $v_1$ satisfies: $L(v_1)=  a_0|\cos s|/2$.
 %
 To  remove secular term in the equation (\ref{rode}) we have to obtain  
 the Fourier expansion for $H(v_0)$ and $v_1$.
 Some computations yield:
 \bqr
  \label{eqabscosF}  |\cos(s)| &= & \frac{2}{\pi} -
      \frac{4}{\pi} \sum_{k=1}^{+\infty} \frac{(-1)^k}{4k^2-1}\cos(2ks),\\
 \label{eqv1F}
    v_1(s) & = & \frac{-a_0}{\pi \omega_0^2}
      \left( 1 - \cos(s)
     -2
  \sum_{k=1}^{+\infty}\frac{(-1)^k}{(4k^2-1)^2}(\cos(2ks) -\cos(s)) 
   \right)
     , 
      \\
\label{eqHF} \nonumber
    H(\cos(s)) & = &  \dis \frac{1}{2} +
      \frac{2}{\pi}  \sum_{k=1}^{+\infty} \frac{(-1)^j}{2j+1}\cos((2j+1)s)
      \label{H(v0)S}.
  \eqr           
 To remove secular term of order one in (\ref{rode}), it suffices to take
 $\alpha_2$ such that:
\bqr \label{evalalpha2}
  0 & = & \dis  \int_0^{2\pi} 
   \left[H(v_0(s))v_1(s)+\alpha_{\bf 2}v_0''(s)+\alpha_1v_1''(s) 
    \right]\centerdot   v_0(s) ds 
 \eqr
For Duffing equation, see \cite{KC68,KC96,Mi}, the source term  involves 
 only few complex exponentials and the calculus of $\alpha_2$ is explicit.  
 For general smooth source term, Fourier coefficients decay very fast. 
 Here, we have an infinite set of frequencies for $v_1$ and $H(v_0)$,
  with only an algebraic  rate of decay for Fourier coefficients.
 So, numerical computations are needed to compute a large number
 of   Fourier coefficients.
  For our first simple example, we can compute explicitly $\alpha_2$.
 After lengthy and  tedious computations involving numerical series, 
 we have from (\ref{evalalpha2}) and 
  (\ref{eqabscosF}), (\ref{eqv1F}),
 to evaluate a 
numerical series
 which yields 
$ \alpha_2=  -3(4\omega_0)^{-2}$ 
thus $\omega_2 = -(2\omega_0)^{-3}$ as we  have already obtained 
 in (\ref{exactomega1}). In more general cases $\alpha_2$ can be computed numerically but not exactly. 
The mathematical result is stated in  Proposition \ref{propexplicit} below. 
The technical proof of the   Proposition \ref{propexplicit}
is postponed to   the appendix.
\\
\\
 We obtain in figure \ref{figspectra}  first modes
 of the 
 Fourier spectra  for 
 $ 
  v_0(\omega_\eps t) + \eps v_1(\omega_\eps t) 
  $
when $a_0=1$.  
%
%
%
\begin{proposition} \label{propexplicit}
 Let $u_\eps$ be  the solution of 
 (\ref{ode}) such that   
   $ 
         u_\eps(0) =a_0 >0 $   
  and
   $ \dot{u}_\eps(0)= 0,
 $ 
 then  there exists $\gamma > 0$, such that,
  for all $t <T_\eps=\gamma\eps^{-1}$, we have the  following expansion
 $ \mbox{ with }
\omega_\eps    =  \omega_0  +\frac{\eps}{4 \omega_0} 
     -\frac{\eps^2}{(2 \omega_0)^3}     $
and $v_1(.) $ is given by  \eqref{eqv1F}:
\bqre 
  u_\eps(t)   & = &a_0 \cos(\omega_\eps t) + \eps v_1( \omega_\eps t)
                       + \mathcal{O}(\eps^2 ) 
 \qquad \mbox{ in } C^2([0,T_\eps],\rit). 
\eqre

\end{proposition}
%
%
%
 Examples from  Proposition  \ref{propexplicit}
  have angular frequency  independent of the amplitude.
  Equation (\ref{ode}) is homogeneous . 
  Indeed, it is a special   case, as we can see in the non homogeneous 
 following cases. In these cases, we assume that the spring is either not in contact with the mass at rest ($b>0$)  or with a prestress at rest ($b<0$).
 \begin{proposition}[Nonlinear dependence of  angular frequency ] \label{PropNdaf} \alali
 Let $b$ be a real number  and let $u_\eps$  be the solution of:
 \bqr \label{1dofgeneral}
    \ddot{u}+ \omega_0^2 u +\eps a(u-b)_+ = 0, &
   &  u_\eps(0)= a_0 + \eps a_1, \;  \dot{u}_\eps(0)=0.
 \eqr
 If $|a_0| > |b|$
then there exists $\gamma > 0$, such that, 
we have the following expansion in $C^2([0,T_\eps],\rit)$
where $T_\eps =  \gamma \eps^{-1}$, 
   $\omega_\eps =  \omega_0 + \eps \omega_1 + \eps^2 \omega_2$
 and $c_k$, $\beta$  are defined by \eqref{ck}, \eqref{beta}
 with $\kappa = b/|a_0|$ :
\bqre 
      u_\eps(t) & = & a_0\cos(\omega_\eps t) +
                    \eps v_1(\omega_\eps t) +
                     \mathcal{O}(\eps^2)
                     \qquad  with  \quad 
     v_1(s)    =   \dis \sum_{k=0}^{+\infty} d_k \cos(k s), 
\\ 
d_k &= &  \frac{- a |a_0|}{\omega_0^2(1-k^2)}
                  c_k  \left(a_0/|a_0|\right)^k,
    \quad \mbox{ for  } 
 k \geq 2, 
\\
      d_0 &= &  -\frac{a |a_0|}{\pi \omega_0^2}
                \left ( \sin(\beta) - \kappa \beta \right), 
    \qquad  
    d_1 =  a_1 -\dis  \sum_{k\neq 1} d_k, \\
    \omega_1 &= & \frac{a }{ 2\pi \omega_0}
           \left (  \frac{\sin (2\beta)}{2} + \beta  
  -  {2\kappa \sin (\beta)} \right ),\\
 \omega_2 &= & \dis \frac{ \omega_1 d_1}{a_0} - \frac{ a}{ \omega_0 \pi a_0} 
  \int_0^{\pi} H(a_0 \cos(s)-b)v_1(s) \cos(s)ds.
\eqre 
 \end{proposition}
Notice that if 
 $|a_0| < |b|$, there is no interaction with the weak unilateral spring.
 Thus the linearized solution is the exact solution.
\\ 
\bpro
   There are two similar cases, $a_0$ positive or negative. 
\\
 \underline{First case}: assume $a_0 > 0$.
 With  the previous notations, the method of strained coordinates yields
 the following equations: 
\bqre 
  v_0''  + v_0 & =&  0, \quad v_0(0)=a_0,\, \dot{v}_0(0)=0
\quad \mbox{ so }  v_0(s) = a_0 \cos(s), 
  \\
-\alpha_0  ( v_1''  + v_1) & = &  a(v_0 - b)_+ -  \alpha_1 v_0
 =  a a_0 (\cos(s) -b/a_0)_+  - \alpha_1 a_0 \cos(s),
 \\
 -\alpha_0  ( r_\eps''  + r_\eps) & = &  
 aH(v_0 - b)v_1 -  \alpha_2 v_0 -\alpha_1 v_1
   + R_\eps.
\eqre
Since  $|\kappa| = | b/a_0 |< 1$, the Fourier coefficient of 
$(\cos(s)-\kappa)_+  = \dis \sum_{k=0}^{+\infty} c_k \cos(k s)$ are:
\bqr 
\label{ck}
c_k= c_k[\kappa] &=&  \frac{1}{\pi}
                \left (\frac{ \sin((k+1)\beta)}{k+1} 
                      +  \frac{ \sin((k-1)\beta)}{k-1}   
               -   \frac{2\kappa \sin(k\beta)}{ k} 
                           \right),
            \; k \geq 2,
\\
  \label{beta}
\beta = \beta[\kappa] & = & \dis \arccos\left(  \kappa \right) \in [0,\pi], 
\\
\nonumber
c_0= c_0[\kappa] &= & \frac{\sin(\beta)  -  \kappa \beta }{\pi },
      \qquad  
c_1= c_1[\kappa] =   \frac{1}{\pi}
           \left ( \frac{\sin(2\beta)}{2} +\beta  -2 \kappa \sin(\beta) \right).
\eqr
The non secular  condition 
 $\dis \int_0^\pi (a(v_0 -b)_+ -\alpha_1 v_0) \cos(s) ds = 0$, 
yields $ \alpha_1 = a\times a_0 \times c_1$.
  Now, we can compute $\omega_1=\alpha_1/(2\omega_0)$
 and the coefficient of the cosinus expansion of $v_1$  
 are  $d_k=  \dis -\frac{aa_0}{\alpha_0} \frac{c_k}{1-k^2}$  for $k \neq 1$.
 The coefficient $d_1$ is then obtained with the initial data $v_1(0)=a_1,$ 
  $\dot{v}_1(0) =0.$ 
$\alpha_2$, is obtained with  the non secular condition for $r_\eps$:
$ 0 = \dis 
\frac{1}{\pi}\int_0^{\pi} (aH(v_0 -b)v_1 -\alpha_2 v_0 - \alpha_1 v_1) \cos(s) ds.
 $
This condition is rewritten as follow
\\
$ \alpha_2 = \dis \frac{2 \omega_0 \omega_1 d_1}{a_0} - \frac{2 a}{\pi a_0} 
  \int_0^{\pi} H(a_0 \cos(s)-b)v_1(s) \cos(s)ds$, 
which gives  $ \omega_2$ since 
 $ \omega_2 =  \frac{\alpha_2 - \omega_1^2}{2 \omega_0}$.
\\
\underline{ Second case}:
 when $a_0 = -|a_0| < 0 $,  by a similar way, we obtain a similar expansion, 
 except that $(v_0(s) -b)_+ = |a_0|(-\cos(s) - \kappa)_+$.
 The Fourier expansion of
 $(-\cos(s)-\kappa)_+ =  \sum_k \tilde{c}_k \cos(ks)$  is simply given 
 by $ \tilde{c}_k = (-1)^k c_k$ since $ - \cos(s)=\cos(s+\pi)$.
 \epro  
  When $|a_0|=|b|$,
 we have  another  asymptotic expansion only valid for shorter time 
 when the unilateral spring slightly interacts with the mass. 
 It is a new feature. 
 \begin{proposition}[Grazing unilateral contact, shorter time  validity]\label{Propcritic}
  \alali
 Let $b$ be a real number, $ b \neq 0$,
  and  consider, the solution $u_\eps$ of problem \eqref{1dofgeneral}.
\\
 If $|a_0|=|b|$  and  $ |a_0+\eps a_1| > |b| $ then  we have  
\\ $ 
  u_\eps(t)  =  (a_0+\eps a_1)\cos(\omega_0 t) + \mathcal{O}(\eps^2),
$
 for  $t \leq \dis T_\eps = \mathcal{O}\left (
 \dis \frac{1}{\sqrt{\eps}} \right)$.
\end{proposition}
Notice that if $ |a_0+\eps a_1| \leq |b| $ 
 then  $u_\eps(t)  =  (a_0+\eps a_1)\cos(\omega_0 t) $ for all time.
\\
The method of strained coordinates gives us 
 the {\it linear }  approximation for $u_\eps(t)$, with $s=t$, i.e.
 $\omega_\eps =1$.
If $|u_\eps(0)| < |b|$, the exact solution is the solution of the linear 
 problem $\ddot{u}+\omega_0^2u=0$.
\\
Otherwise,  if $|u_\eps(0)| > |b|$,
since, $|b|$ is the maximum of $v_0(s)=a_0\cos(s)$, a new phenomenon appears, 
 during each period, $|u_\eps(t)| > |b|$ on  interval of time 
 of order $\sqrt{\eps}$ instead of $\eps$. 
Then   $T_\eps$ is smaller  than in Proposition \ref{propexplicit}.
\\
  To explain this  phenomenon, we give precise estimates of the remainder 
 when we expand $(v_0+\eps v_1 + \eps^2 r_\eps)_+$ in the next section,
   see  Lemmas \ref{LDLL} and \ref{LHzi} below.

    \section{Expansion of $(u+\eps v)_+$ \label{sexp+}}  

 We give some useful lemmas to perform asymptotic expansions and 
to estimate precisely the remainder  for the   piecewise linear map 
    $ u \rightarrow u_+= \max(0,u).$ 
  
 \begin{lemma}
 {\bf [Asymptotic expansion  for $(u+\eps v)_+$
  ] } 
\label{LDLL}
Let be $T>0$, $u, v$ two real valued functions defined on $I=[0,T]$,and 
  $H$ be the Heaviside  step function then 
\bqr  \label{equa:DLL}
     (u+\eps v)_+  & = & (u)_+ + \eps H(u)v 
    + \eps \chi_\eps(u,v) ,
\eqr
 where  $\chi_\eps(u,v)$ is a non negative piecewise linear 
    function and  1-Lipschitz with respect to $v$.
\\ 
Let be $M>0$,
 $J_\eps =\{t \in I, |u(t)|\leq \eps M \},$
$\mu_\eps(T)$ the measure of the set $J_\eps$. 
\\
 If  $ |v(t)|\leq M $  for any $t \in I$ then  
  \bqr \label{eq:DLL} 
        |\chi_\eps(u,v)| \leq   |v| \leq M, &  
           & 
       \dis \int_0^T \left |\chi_\eps(u(t),v(t))  \right|dt   \leq 
                     M \mu_\eps(T).
  \eqr 
\end{lemma}
The point in inequality (\ref{eq:DLL}) is  the remainder $\eps \chi_\eps$ is only 
 of order $\eps$ in $L^\infty$ but of order $\eps \mu_\eps$ in $L^1$. 
 In general, $\mu_\eps$ is not better than a constant, 
take for instance $u\equiv 0$. Fortunately, 
 it is proved below that  $\mu_\eps $ is often of order 
 $\eps$, and for some critical cases of order $\sqrt{\eps}$.
\\
\bpro Equality (\ref{equa:DLL}) 
defines $\chi_\eps$ and can be rewritten as follow: 
  \bqr  \label{eqchi}
       \chi_\eps(u,v) 
 &=&
     \dis  \frac{  (u+\eps v)_+  -    u_+ - \eps H(u)v}{\eps}.
\eqr   %
  So,  $\chi_\eps$ is non negative since 
  $u \ra u_+$ is a convex function.
We also  easily see that the map 
$(u,v) \rightarrow \chi_\eps(u,v) $ is piecewise linear, continuous 
 except on the line $u=0$ where $\chi_\eps $ has a jump  $-v$.
 This jump comes from the Heaviside step function. 
An explicit computations gives us the simple and useful formula: 
   $   0 \leq \eps  \chi_\eps(u,v) 
      = |u+\eps v|  H(|\eps v|- |u+\eps v|)   $. 
  We then have immediately $ 0 \leq \chi_\eps(u,v)  \leq |v| $.
  Let $u $ be  fixed, then $v \rightarrow \chi_\eps(u,v)$
 is one Lipschitz with respect 
  to $v$.
Furthermore,
the support of $\chi_\eps $ is included in $J_\eps$,
    which  concludes the proof. 
 \epro  
\noindent 
Now we  investigate the size of  $\mu_\eps(T)$, see \cite{BJ,J4} 
for similar results about $\mu_\eps(T)$ and other applications.
With notations
 from  Lemma \ref{LDLL} we have. 
\begin{lemma}[Order of $\mu_\eps(T)$]\label{LHzi}
 Let $u$ be  a smooth  periodic function, $M$ be a positive constant   and 
  $\mu_\eps(T)$ the measure of the set
 $J_\eps= \{t \in I, |u(t)|\leq \eps M \}$.
\\ 
    If $u$   has only  simple roots on $I=[0,T]$
then  for some  positive $C$,
  $ 
 \mu_\eps(T)  \leq C\eps \times T .
 $
\\
More generally, if $u$  
    has also double   roots 
then  $  \mu_\eps(T) \leq C\sqrt{\eps} \times T.$ 
\end{lemma}
The measure of such set $J_\eps$ implies 
many applications in averaging lemmas, 
for a characterization of $\mu_\eps$ in a multidimensional framework
 see \cite{BJ,J4}.
\\
Notice that
   any non zero  solution $u(t)$ 
  of any linear homogeneous second order ordinary differential equation
  has  always simple  zeros, thus 
 for any constant $c$ the map $t \rightarrow u(t)-c$ has at most double roots.
\medskip \\
\bpro 
 First assume  $u$ only has simple roots on a  period $[0,P]$, 
 and let $Z = \{t_0 \in [0,P],u(t_0)=~0\}$. 
 The set $Z$ is  discrete  
  since $u$ has only  simple roots which implies that 
     roots of $u$ are isolated.
 Thus $Z$ is a finite subset of $[0,P]$: 
   $Z=\{t_1, t_2, \cdots,t_N \}$.
We can choose an open  neighborhood $V_j$ of each  $t_j$
 such that $u$ is a  diffeomorphism on $V_j$
  with derivative $|\dot{u}| > |\dot{u}(t_j)|/2$.
On  the compact set $K=[0,P]-\cup V_j$,
  $u$ never vanishes, then
   $\dis \min_{t \in K} |u(t)|=\eps_0 >0.$
 Thus, we have for all $\eps M < \eps_0$, the length of $J_\eps$ in $V_j$
 is  
   $\dis  |V_j \cap J_\eps|  \leq \frac{4 \eps M}{|\dot{u}(t_j)|} $.
  As $\mu_\eps$ is additive ($\mu_\eps(P+t)= \mu_\eps(P)+\mu_\eps(t)$),
  its growth is linear. 
 Thus, for the case with simple roots, we get 
$\mu_\eps(T) = \mathcal{O}(\eps T)$.  
\\
For the general case with double roots, 
on each small neighborhood of $t_j$: $V_j$, 
we have with a Taylor expansion,
  $|u(t_j+s)| \geq  d_j|s|^l$, with $ 1 \leq l \leq 2$, $d_j > 0$, so,
 $\dis  |V_j \cap J_\eps| \leq 2 (\eps M/d_j)^{1/l}  $, then 
$\dis \mu_\eps(P) =  \mathcal{O}(\sqrt{\eps}), $
 which  is enough
to  conclude the proof.
\epro

    \section{Several  degrees  of freedom \label{sNdof}}   
 Now, we investigate the case with $N$ masses.
 We use, the method of strained coordinates in three cases.
 We present  the formal computations for each expansion. The mathematical proofs 
 are postponed in the Appendix.
\\ 
 In subsection \ref{ssNdof1}, the initial condition is near an eigenvector
 such that 
 the approximate solution  stays periodic.
  We give such  initial condition near an eigenvector
 in subsection \ref{sec:annm} to get an approximate nonlinear normal mode
 up to the order $\eps^2$.
 Finally,  in subsection \ref{ssNdof2}, all modes are excited. 
 An extension of the method of strained coordinates is still possible but 
 only at the first order with less accuracy. 
\\
The system studied is the following:

$M\ddot{\tilde{U}} + K \tilde{U}  +\eps (\tilde{A}\tilde{U} -B)_+=0$, 
where, for each component, 

   $[(\tilde{A}U-B)_+]_k = \dis\left( \sum_{j=1}^N \tilde{A}_{kj}u_j -b_k\right)_+$,
$M$ is  a  $N \times N$ 
mass matrix,
$K$ is  the stiffness  matrix they are both  symmetric definite positive.  
$\tilde{A}$ and $B$ are matrices which involve the rigidity of unilateral springs and their position with respect to the masses.
   For  such  a system, endowed with a natural convex 
 energy for the linearized part, 
  we can control the $\eps$-Lipschitz nonlinear term for $\eps$ small enough up to large time.
 So for $\eps <<1$ the solutions remain bounded for  time of the order $\eps^{-1} $.
\\
We introduce the matrix $\Phi$ of generalized eigenvectors:
$ K \Phi = M\Phi \Lambda^2 $ with $\Lambda$ positive diagonal matrix of eigenvalues, 
 $\Phi^T M \Phi = Id$,  and set $ \tilde{ U}=\Phi U$, $A=\tilde{A} \Phi$,  the system may be written:

\begin{eqnarray}\label{sysMK+}
   \ddot{U} + \Lambda^2 U  &=& -\eps \Phi^T (AU -B)_+ .
 \end{eqnarray}

\subsection{Initial condition near  an eigenvector, 
     \label{ssNdof1}}   
  For the system (\ref{sysMK+}), we take  
 an initial condition near  an eigenmode of the linearized system
denoted for instance by  index $1$  .
\bqr    \label{initialNNM}
  \left \{ 
  \begin{array}{ccccc} u_1^\eps(0) & =& a_0+\eps a_1, &\dot{u}_1^\eps(0)=0, & \\
   u_k^\eps(0)&=&0\;\; + \eps a_k,  &\dot{u}_k^\eps(0)=0,&   for \; k \neq 1. 
 \end{array}
   \right.
\eqr 
We impose $a_2,\cdots, a_N$ later  to have a periodic approximation, but 
 $a_1$ is a free constant as $a_0$. 
It is a key point to apply the method of strained coordinates.
\\
 We use  the same time $ s= \omega_\eps t$  for each component 
  and  the following  notations. 
\begin{eqnarray*}
\begin{array}{cclccl}
   \omega_\eps   & =&  \omega_0 + \eps \omega_1+\eps^2 \omega_2,
  &
   \omega_0 & = & \lambda_1,
  \\ 
    (\omega_\eps)^2 & = & \alpha_0 + \eps \alpha_1+\eps^2 \alpha_2 
 + \mbox{\rm O}(\eps^3), 
  &
       \alpha_0 &=& \omega_0^2= \lambda_1^2,
  \\
          \alpha_1&=& 2\omega_0 \omega_1,
    &  \alpha_2&=& \omega_1^2 + 2 \omega_0 \omega_2,
   \\
   u_j^\eps ( t)   & = &   v_j^\eps (s) 
   =    v_j^0(s) +\eps v_j^1(s)+ \eps^2  r_j^\eps(s),
 & j&=&1,\cdots,N . 
 \end{array}
 \end{eqnarray*}
Replacing, this ansatz in the System (\ref{sysMK+}) we have in variable $s$,
\bqre \label{dasNLNM}
  ( \omega_\eps)^2 (v_k^\eps)''  + \lambda_k^2 v_k^\eps  
& =&
  -\eps \dis \sum_{l=1}^N \Phi_{lk}  \left( \sum_{j=1}^N A_{lj}
v_j^\eps\left(s \right) -b_j\right)_+,
 \eqre
 and then performing the expansion for all $k \in\{1,\cdots,N\}$, 
\bqre  \label{eqNLM0} 
  L_k v_k^0  &=& \alpha_0 (v_k^0)''  + \lambda_k^2 v_k^0  
 = 0 ,
  \\
\label{eqNLM1}
-L_k v_k^1  
 &= &  \dis  \sum_{l=1}^N \Phi_{lk}  \left(\sum_{j=1}^N A_{lj}
v_j^0 -b_l\right)_+ + \alpha_1(v_k^0)'', \\ 
\label{eqNLM2}
 -L_k r_k^\eps  
 & =&
     \dis \sum_{l=1}^N \Phi_{lk}  H\left( \sum_{j=1}^N A_{lj}
v_j^0 -b_l\right)
\left( \sum_{j=1}^N A_{lj}
v_j^1\right)
+ \alpha_2 (v_k^0)''+ \alpha_1 (v_k^1)''
+R_k^\eps.
\eqre
  First we have $v_1^0(s)= a_0 \cos(s)$.
\\
  Equations for $v_k^0$, for all $k\neq 1$, with zero initial data 
 give us $ v_k^0 = 0$.
\\ 
 In equation for $v_1^1$,
we remove the secular term in the  right hand side,
\bqre 
-\alpha_0( (v_1^1)''+v_1^1)  
 =
  \dis \sum_{l=1}^N \Phi_{l1}  \left(  A_{l1}
v_1^0 -b_j\right)_+ + \alpha_1(v_1^0)''=r.h.s.
 &
 \quad 
 & v_1^1(0)=a_1, \; (v_1^1)'(0)=0 .  
\eqre
The orthogonality of the $r.h.s$ with $\cos(s)$  defines $\alpha_1$ with \eqref{eqalpha1def}.
For instance, if $b_1  =0$ and $A_{11}>0$, 
 we have as in Proposition \ref{PropNdaf}, 
 $ 2\alpha_1= \dis \sum_l \Phi_{l1} A_{l1}$ 
 and $ \omega_1=\dis \frac{\alpha_{1}}{2\lambda_1}$.
\\
Now, $\alpha_1$ is fixed, 
so $v_1^1$ is a well defined  even $2\pi$ periodic function.
\\
 Then, for $k \neq 1$, $v_k^1$  is the unique $2\pi$  periodic  solution 
of  the  simplified equation,
\bqr \label{eqN1} \dis 
 - L_k v_k^1  
   &=&   \dis \sum_{l=1}^N \Phi_{lk}  \left(  A_{l1}
v_1^0 -b_l \right)_+.
\eqr
Such a function exists and is unique  if  $\lambda_k \nin \lambda_1 \zit$.
Furthermore $v_k^1$ is an even function as the right hand side of  equation 
(\ref{eqN1}).  Then $a_k$ is given by $v_k^1(0)$ and $(v_k^1)'(0)=0$ for all $k\neq 1$. \\
   The term $r_1^\eps$, with null initial data,
 has  a simplified equation 
 since $v_k^0 \equiv 0$ for all $k\neq 1$,
 \bqre 
 -L_1 r_1^\eps  
 & =&
     \dis \sum_{l=1}^N \Phi_{l1}  H\left(  A_{l1}
v_1^0 -b_l\right)
\left(   \sum_{j=1}^N A_{lj}
v_j^1\right)
+ \alpha_2 (v_1^0)''+ \alpha_1 (v_1^1)''
+R_1^\eps.
 \eqre 
 Now we  can compute numerically 
 $\alpha_2$ to avoid secular term in the right hand side, $R_1^\eps$ excepted,
 with the following condition, 
\bqre 
 0 &= & 
   \int_0^\pi \left[ \dis \sum_{l=1}^N \Phi_{l1}  H\left(  A_{l1}
v_1^0 -b_l\right)
\left( 
\sum_{l=1}^N A_{l1}
v_l^1\right)
+ \alpha_2 (v_1^0)''+ \alpha_1 (v_1^1)''
 \right] \cdot \cos(s) ds.
 \eqre
Rewriting this condition, we obtain an equation for $\alpha_2$  
in Theorem \ref{Th41}  below.
 \\
  For each $k\neq 1$, $\lambda_k \nin\lambda_1\zit$, 
   so  $r_k^\eps$ stays bounded 
  for large time. Indeed there is no resonance of the order one at the first 
 order in equation satisfied by  $r_k^\eps$. 
 This  is the  technical  part of the proof 
   to validate rigorously
  and  to find the time of validity of 
  such  asymptotic expansion.
The complete proof to bound $(r_1^\eps, \cdots, r_N^\eps)$ for large time
 is to be found in the Appendix, subsection \ref{ssbr}.
 Now we state our result with previous notations. 

\begin{theorem}\label{Th41}
The Lindstedt-Poincar\'e expansion is valid
 on $(0,T_\eps)$, with $T_\eps \rightarrow  + \infty$ when $\eps \rightarrow 0$
 under assumption $ \{\lambda_2,\cdots,\lambda_N\} \cap \lambda_1 \zit =\emptyset$:
 \bqre
\label{eq:uup=0}
  \left\{ 
   \begin{array}{cccccccc}
   u_1^\eps(t) &  = &  v_1^0(\omega_\eps t) & +& \eps v_1^1(\omega_\eps t)
                   &+& \mathcal{O}(\eps^{2}),& 
                \\
 u_k^\eps(t) &  = & 0 & +& \eps v_k^1(\omega_\eps t)
                   &+& \mathcal{O}(\eps^{2}),&   k \neq 1, 
   \end{array}
    \right.
 \eqre
where $v_1^0(.)$, $\alpha_1$, $\omega_1$, $v_1^1(.)$,
  $v_k^1(.)$ and  $a_k$ for  $k \neq 1$,  $\alpha_2$, $\omega_2$
 are {\rm successively} defined  as follows: 
\bqre 
v_1^0(s) &=& a_0 \cos(s), \\
  \alpha_1 &= & \dis \frac{2}{a_0\pi}
\int_0^\pi \sum \Phi_{l1}
      (A_{l1}v_1^0(s) - b_1)_+ \cos (s) ds, \label{eqalpha1def} 
  \mbox{ then } \omega_1 =  \frac{\alpha_1}{2 \omega_0}, \\ 
-L_1 v_1^1   & = &  \dis\left(  A_{11}
v_1^0 -b_1\right)_+ + \alpha_1(v_1^0)'',
\quad  v_1^1(0)=a_1, \; (v_1^1)'(0)=0,
\\
 v_k^1 & \mbox{ be } &
 \mbox{the unique } 2\pi \mbox{ periodic solution  of } \eqref{eqN1}  \mbox{ and } a_k:=v_k^1(0),  \mbox{ for } k \neq 1, 
\\  
 \alpha_2 &=&   \left( \frac{2}{a_0\pi}  \sum_{l=1}^N  \Phi_{l1} 
\int_0^\pi  H_{l1}(s)\cos (s) ds \right) 
 + \alpha_1 \int_0^\pi (v_1^1)'' cos(s) ds,
\\
 &   \mbox{with} & H_{l1}(s) =    H(A_{l1}v_1^0(s) - b_l)
     \left(A_{l1}v_j^1(s)+ \sum_{k\neq 1}^N A_{1k}v_k^1(s) \right),   
\\
  \omega_{\eps} & =&  \omega_{0}+ \eps \omega_{1} + \eps^{2}\omega_{2},
\quad \mbox{ where }
 \omega_2 =  \frac{\alpha_2 - \omega_1^2}{2 \omega_0} \qquad 
\eqre
and
 $L_k$ be the differential operator $\dis \lambda_1^2 \frac{d^2}{ds^2} + \lambda_k^2 $.
\\
 Furthermore, if $(A_{j1} v_1^0-b_j)$ has got only simple roots 
for all $j=1,\cdots,N$, 
\\ then   $T_\eps = O(\eps^{-1} )$,
\quad else  $T_\eps= O(\eps^{-1/2} )$.
\end{theorem}
In the theorem,  $v_k^1$ is classically obtained by  a Fourier series. 
We give some indications of its initial condition in the 
 next subsection \ref{sec:annm}.

\subsection{Approximate non linear normal mode}
\label{sec:annm}

The special initial conditions  of the previous subsection
  can   be explicited
 in order to find a solution where all the components are in phase at the
same frequency.
Indeed we shall obtain an approximate curve of initial conditions for  which the solution is periodic up to the order $\eps$ for a time of the order $\eps^{-1}$ or $\eps^{-1/2}$: this is up to the second order approximation a non linear normal mode in the sense of Rosenberg \cite{Rosenberg}; see \cite{JPS} for a computation of non linear normal mode with unilateral contact and \cite{nnm-kpgv} for a synthesis on non linear normal modes.

\begin{corollary}[Explicit initial condtion for the approximate NNM] \alali
\label{thm:app-nnm}
Let  $a_0 \neq 0$  be fixed, $k \neq 1$,  
$A_{j1}\neq 0$, $\kappa_j=\frac{b_j}{A_{j1}a_0} $,
and $c_l[\kappa]$ defined by   formula \eqref{ck}. 
 $a_k$ from \eqref{initialNNM}   are  computed explicitly in the following cases:
\begin{enumerate}
\item if $b_j=0$ for $j=1,\cdots,N$, then 
  \begin{equation}
    \label{eq:akpourbk0}
    a_k= \sum_{j=1}^N \Phi_{jk} \left( 
                    \frac{| A_{j1}a_0|}{2(\lambda_1^2-\lambda_k^2)}
 - \frac{|A_{j1}a_0|}{\lambda_k^2 \pi} 
+\frac{ | A_{j1}a_0|}{\pi} \sum_{l=1}^{+\infty} \frac{(-1)^l}{(4 l^2\lambda_1^2 -\lambda_k^2)(4l^2 -1)}
 \right)
  \end{equation}

\item 
   if  $0 < b_j $,
      $|\kappa_{j}|<1$, $a_0 \; A_{j1}<0$ for all $ j$
  then 
  \begin{equation}
    \label{eq:akpourbk0bkak1a0pp1p}
    a_k=  -\sum_{j=1}^N   \Phi_{jk}  |a_0A_{j1}| 
\left [ 
\sum_{l=1}^{+\infty} \frac{ (-1)^lc_{l}[-\kappa_j]}{l^2 \lambda_1^2 -\lambda_k^2}
\right]
\end{equation}
\item
 if $0 < |\kappa_j|<1$, and $a_0 \; A_{j1}>0$ for all $j $ then 
  \begin{equation}
    \label{eq:akpourbk0bkak1a0pp1m}
    a_k= \sum_{j=1}^N \Phi_{jk} a_0A_{j1} \left [ 
\sum_{l=1}^{+\infty} \frac{c_{l}[\kappa_j]}{l^2 \lambda_1^2 -\lambda_k^2}
\right]
  \end{equation}
\end{enumerate}

\end{corollary}
Thanks to Theorem \ref{Th41}, such initial data given by $(a_k)_k$ determine the approximate NNM.
 Notice that there is no condition on $a_1$.
The other numerous  cases 
  may solved similarly.  
\\
\\
\bpro
The principle of the proof is simple: $v_k^1$
 is  the periodic solution of the differential equation  \eqref{eqN1} 
and $a_k=v_k^1(0)$ has  to be determined in order that the function $v_k^1$  has an angular frequency equal to one. Solution of \eqref{eqN1} is 
$
  \label{eq:v_k^1=b...}
  v_k^1=A \cos\left(\frac{\lambda_k}{\lambda_1} s\right)
  +B \sin\left(\frac{\lambda_k}{\lambda_1} s\right) + w_k^1(s),
$
where $w_k^1$ is a particular solution associated to the right hand side which is of angular frequency equal to $1$. Note that $B=0$ as the initial velocity is null. We can get a function of angular frequency equal to $1$ by setting $a_k=w_k^1(0)$. This condition may be written explicitly with formulas 
\eqref{ck} which provides the expansion in Fourier series; formulas 
\eqref{eq:akpourbk0}, \eqref{eq:akpourbk0bkak1a0pp1p},     \eqref{eq:akpourbk0bkak1a0pp1m}   are then derived easily successively.
\begin{enumerate}
\item If $b_k=0$, $ k\ne 1$, \eqref{eqN1} is written:
  $
    -L_k v_k^1= \sum_{j=1}^N \Phi_{jk} \left(
A_{j1} a_0 \frac{\cos(s)}{2} +  | A_{j1} a_0 |  \frac{|\cos(s)|}{2} \right)
 .$
We use formula \eqref{eqabscosF} to get the particular solution
$w_k^1$ and then formula \eqref{eq:akpourbk0}: 
\\ $
  w_k^1(s) =  \sum_{j=1}^N\Phi_{jk} \left(
A_{l1} a_0 \frac{\cos(s)}{2(\lambda_1^2 -\lambda_k^2) }-
\frac{|A_{j1} a_0|}{\lambda_k^2 \pi} +
\frac{2|A_{j1} a_0|}{\pi} \sum_{l=1}^{+\infty} \frac{(-1)^l}{4l^2 \lambda_1^2-\lambda_k^2 } \frac{\cos(2l s)}{4l^2 -1} \right).
$
\item For the second case, 
\eqref{eqN1} is written:
$
      -L_k v_k^1= -\sum_{j=1}^N \Phi_{jk}
A_{j1} a_0 \left (-\cos(s) + \kappa_j  \right )_+.  
$
 We use \eqref{ck} to obtain 
$w_k^1(s)=  -\sum_{j=1}^N \Phi_{jk}
|A_{j1} a_0 | \left [ 
\sum_{l=1}^{+\infty} \frac{c_{jl} \cos(ls)}{l^2 \lambda_1^2 -\lambda_k^2} \right ]
$
where $ c_{jl}= c_l[-\kappa_j]  $.
\item 
 For the third case \eqref{eqN1} is written:
$
       -L_k v_k^1= \sum_{j=1}^N\Phi_{jk}
A_{j1} a_0 \left (\cos(s) - \kappa_j \right )_+
 $ 
  from which 
$ w_j^1(s)= \sum_{j=1}^N\Phi_{jk}
a_0 A_{j1} \sum_{l=1}^{+\infty} \frac{c_{jl} \cos(ls)}{l^2 \lambda_1^2 -\lambda_k^2} $
 where  $c_{jl}=c_l[\kappa_j] $.
\end{enumerate}

\epro

\subsection{Numerical results of NNM}

\subsubsection{Using numerically Lindstedt-Poincar\'e  expansions}
Here we use the previous results and  compute numerically a solution of system 
\eqref{sysMK+} using the approximation \eqref{eq:uup=0}:
$u^{\eps}(t)=v^0(\omega_{\eps} t) +\eps v^1(\omega_{\eps} t) +O(\eps^2)$ 
with the initial conditions of theorem \ref{thm:app-nnm}. The first term $v^0$ is easy to obtain; for the second term $v^1$ an explicit  formula is in principle possible using Fourier series such as for one degree of freedom but it is cumbersome so we choose to compute $v^1$ by solving numerically 
 \eqref{eqN1}
 with a step by step algorithm; we use as a black-box the routine ODE of SCILAB \cite{scilab} to solve equations of theorem \ref{Th41} after  computing by numerical integration $\alpha_1$. We show numerical results for a system of the type:
\begin{gather} \label{sysMKF}
  M \ddot{X} + K X + \eps F(X)=0
\end{gather}
we still denote $\lambda_j^2$ the eigenvalues and $\phi_j$ the eigenvectors of the usual generalized eigenvalue problem
 $ 
  K \phi_j -\lambda_j^2  M \phi_j=0 \quad \text{with} 
^t \phi_k M \phi_j=\delta_{kj}.
$
 We set:
$ 
  X  = \sum_j u_j \phi_j =\underline{ \underline{\phi}} \underline{u}.
$
In this basis, the system may be written componentwise:
$
\ddot{u_k}  +\lambda_k^2  u_k + ^t \phi_k  \eps F(\underline{ \underline{\phi}} \underline{u})  =  0.
$
\\
We illustrate a  simple local non linearity in the system \eqref{sysMKF},
 with the following nonlinearity 
\bqre
  F(X) & = &  (X_1-\beta_1) M \phi_1 
    =\left(\sum_j u_j \phi_{j1}  -\beta_1 \right)_+  M \phi_1.
\eqre
The system \eqref{sysMKF},  written in the basis of the eigenvectors, 
simply becomes:
\bqre 
  \ddot{u_1}  +\lambda_1^2 u_1 
 + \eps \left(\sum_j\phi_{j1} u_j - \beta_1 \right)_+   = 0   , & \mbox{ and for } k \neq 1   & 
 \ddot{u_k}  +\lambda_k^2 u_k = 0.
\eqre

%
We find in figure \ref{fig:lp-ener0_03} a numerical example of the Lintsted-Poincar\'e approximation for 5 degrees of freedom with $\eps=0.063$ and with an energy of $0.03002$. The left figure shows the 5 components of the solution with respect to time; the right figure, the solution in the configuration space: abscissa component 1 and ordinate components 2 to 5; these lines are rectilinear like in the linear case but the non symmetry may be particularly  noticed on the smallest component which corresponds to the mode where the non linearity is active. 
\begin{figure}[!ht]
  \centering
\begin{center}
  \includegraphics[ width=7cm, height=7cm,angle=-90]{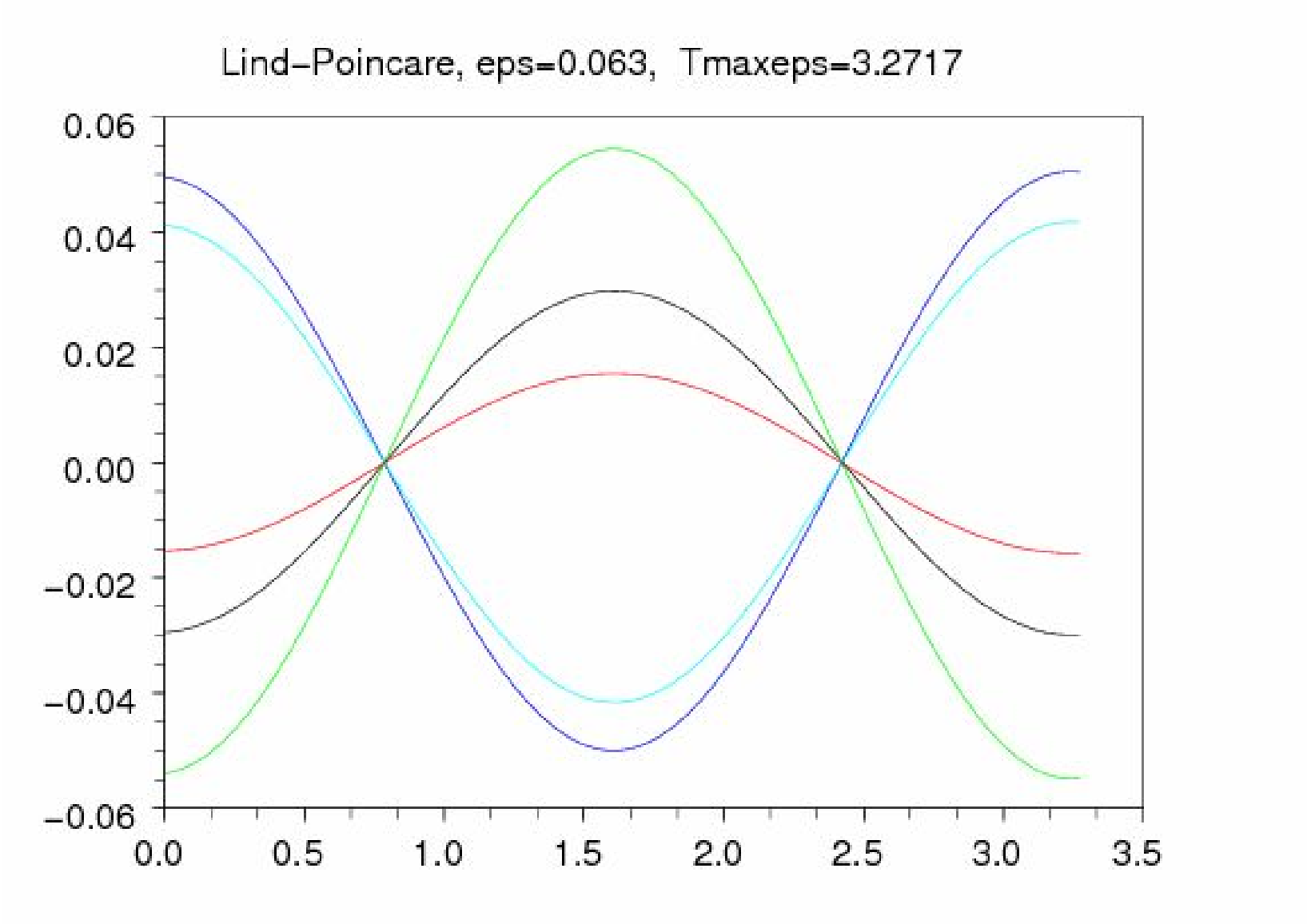}
  \includegraphics[ width=7cm, height=7cm, angle=-90]{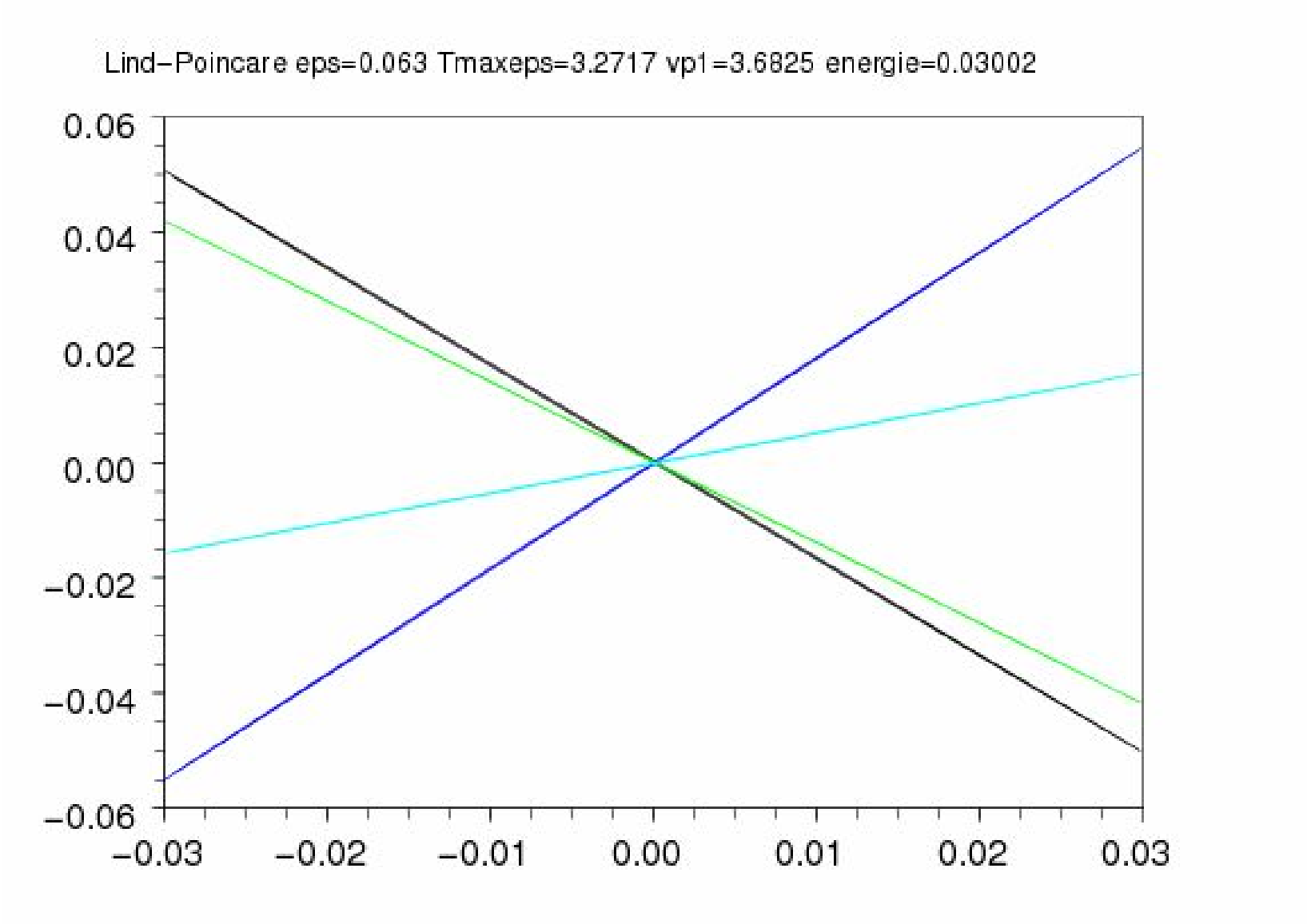}
\caption{Lindstedt-Poincar\'e, energy=0.03, 5 dof; left: components with respect to time; right: in configuration space  \label{fig:lp-ener0_03} }
\end{center}
\end{figure}

\subsubsection{Using optimization routines} \label{sssubrotines}
We also find   in figure  4  
a numerical example with the same energy of $0.03002$; it is computed with a purely numerical method described below. We notice that the solution is quite similar in both cases.

\label{sec:opti-routi}
\begin{figure} \label{fig:opt-ener0_03}
  \centering
\begin{center}
  \includegraphics[ width=7cm, height=7cm, angle=-90]{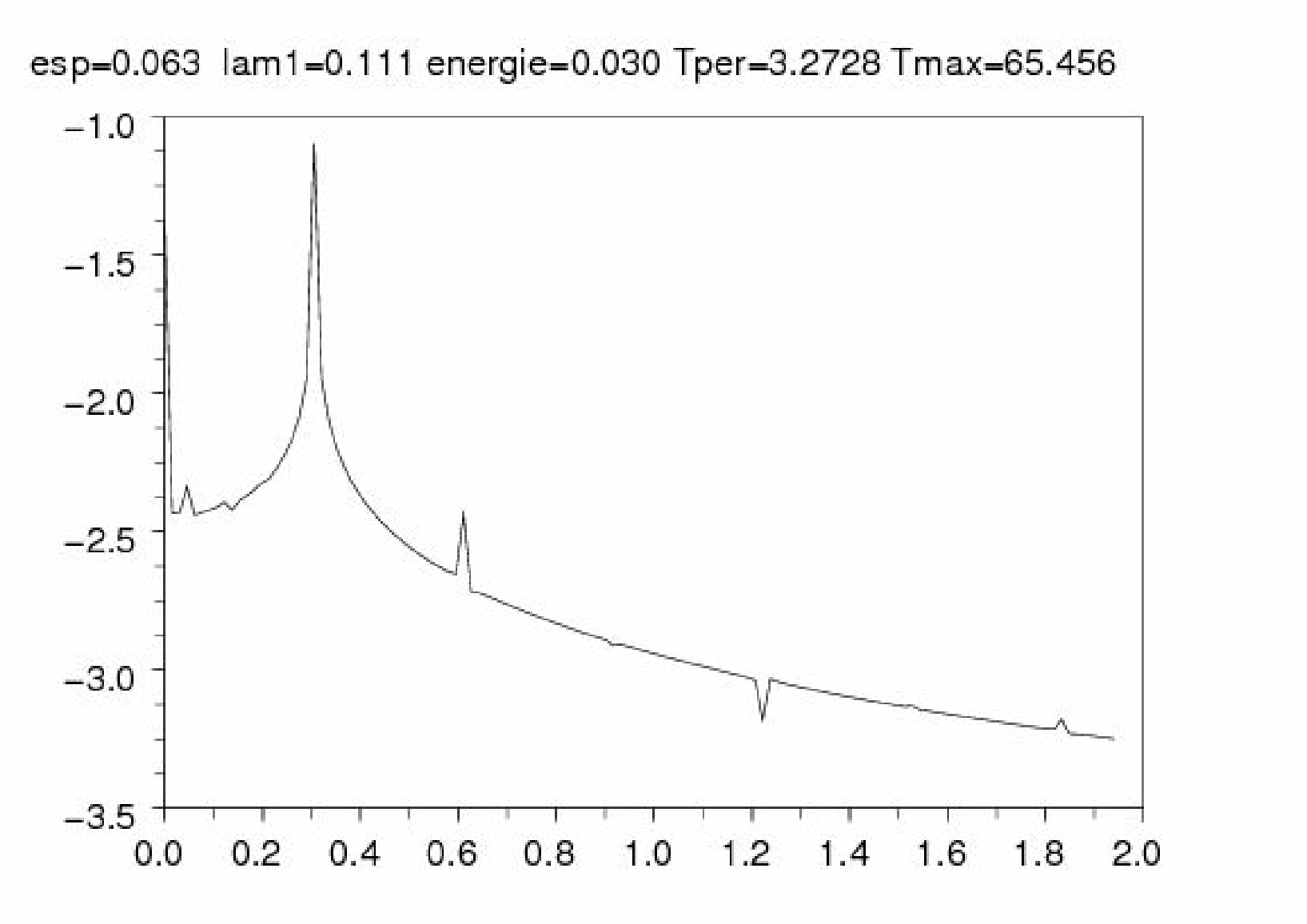}
  \includegraphics[ width=7cm, height=7cm,angle=-90]{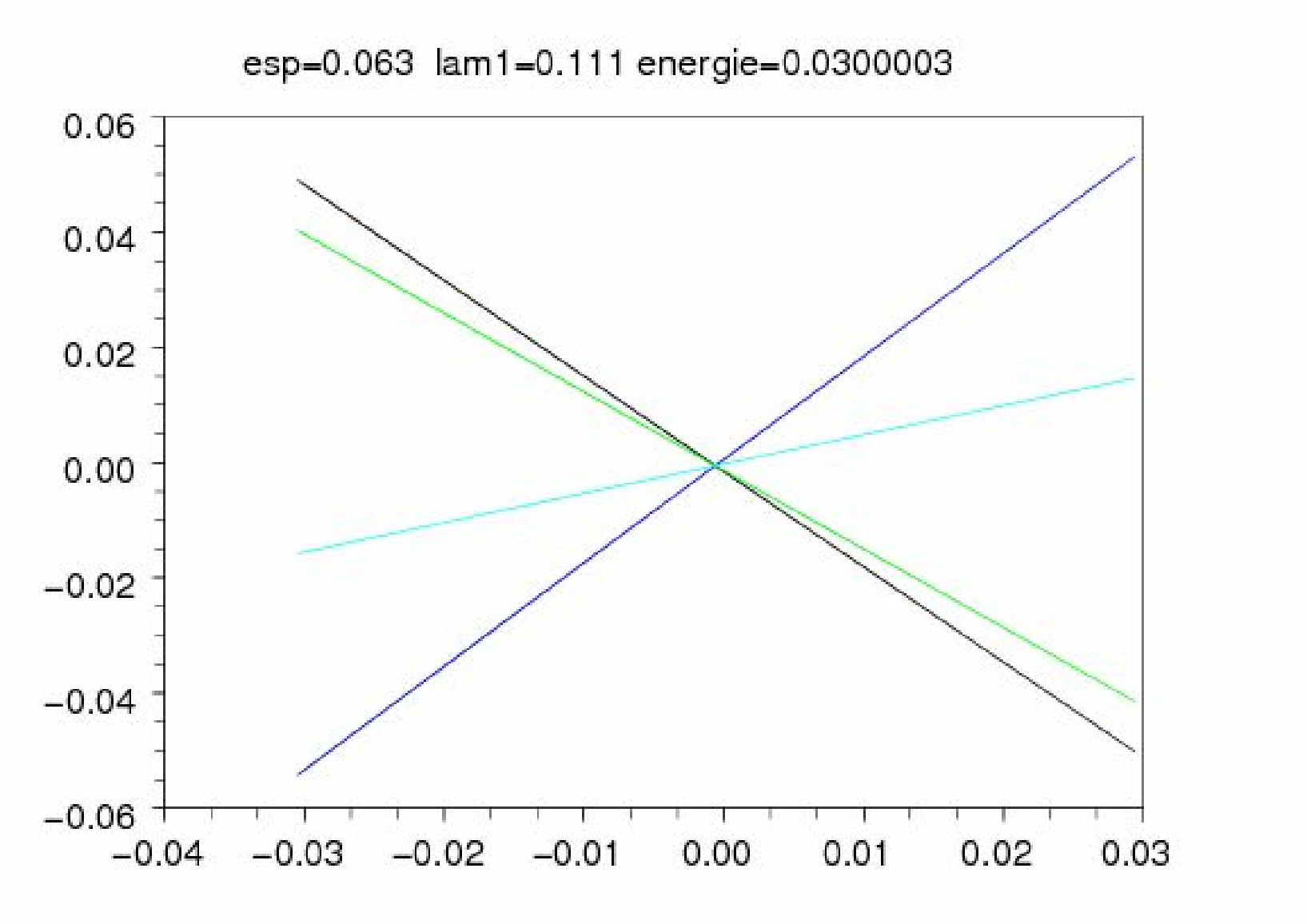}
\caption{Continuation and Powell hybrid, energy=0.03, 5 dof; left:Fourier transform ; right: in configuration space  \label{}}
\end{center}
\end{figure}

The numerical expansions of the previous subsection gives valid results for $\eps$ small enough; in many practical cases such as  \cite{HR}, $\eps$ may be quite large; in this case, it is natural to try to solve numerically the following equations with respect to the period $T$ and the initial condition $X(0)$.
\bqre 
\label{eq:X0=XT}
  X(0)=X(T), \quad \dot{X}(0)=\dot{X}(T) \quad
E(X)=e
\eqre 
In other words, we look for a periodic solution of prescribed energy; this last condition is to ensure to obtain an isolated local solution: the previous expansions show that in general, the period of the solution depends on its amplitude prescribed here by its energy. To try to solve these equations with a black-box routine for nonlinear equations such as ``fsolve'' routine of SCILAB \cite{scilab}(an implementation of a modification of Powell hybrid method which goes back to \cite{powell64}) in general fails to converge.
Even in case of convergence, we should address the question of link of this solution with normal modes of the linearized system.
\par So we prescribe that $e=c\eps$ and for $\eps\rightarrow 0$, the solution is tangent to a linear eigenmode. In the case where all the eigenvalues of the linear system are simple, we define $N$ (the number of degrees of freedom) non linear normal modes for which, it is reasonable to conjecture that they correspond to isolated solutions of  \eqref{eq:X0=XT} at least for small $\eps$ if we enforce for example $\dot{X}(0)=0$.
\paragraph{Algorithm}
This definition of the solution of \eqref{eq:X0=XT} tangent to a prescribed linear eigenmode provides a simple way of numerical approximation: using a continuation method coupled with a routine for solving a system of non linear equations.
Define:
\bqre
\mathcal{F}(\eps,X_0,X_1, T) & = &[X(T)-X_0,\dot{X}(T)-X_1, E(X)-c\eps], 
\eqre
where X is a numerical solution of the differential system
 $  \dis \left \{ \begin{array}{c} M \ddot{X} + K X + \eps F(X)=0 \\ 
                 X(0)=X_0, \quad \dot{X}(0)=X_1 \end{array}
\right.$
\begin{verse}
choose a small initial value of $ \eps$ and an increment $ \delta$\\
choose an eigenvector $\phi_j$ \\
$X_0(0)= A_{\eps} \phi_j$, $X_1(0)=  B_{\eps} \lambda_j \phi_j$\\
with $E(X_0(0),X_1(0))=c \eps$ \\
\end{verse}
\begin{verse}
{\bf for}  iter=1:itermax
$$\eps=\eps + \delta$$
 with $ (X_0(iter-1),X_1(iter-1)) $
as a first approximation, solve for $(X_0(iter),X_1(iter))$,
\bqre
 \mathcal{F}(\eps,X_0,X_1,T)& =& 0  
\eqre
{\bf if } $||\mathcal{F}(\eps,X_0,X_1,T)||> tolerance$ {\bf then}
 $\eps=\eps - \delta, \quad\delta=\delta/2 $\\
$\quad$   {\bf endif} \\
{\bf endfor}
\end{verse}

\begin{figure}[!ht] 
\begin{center}
\includegraphics[ width=7cm, height=6cm]{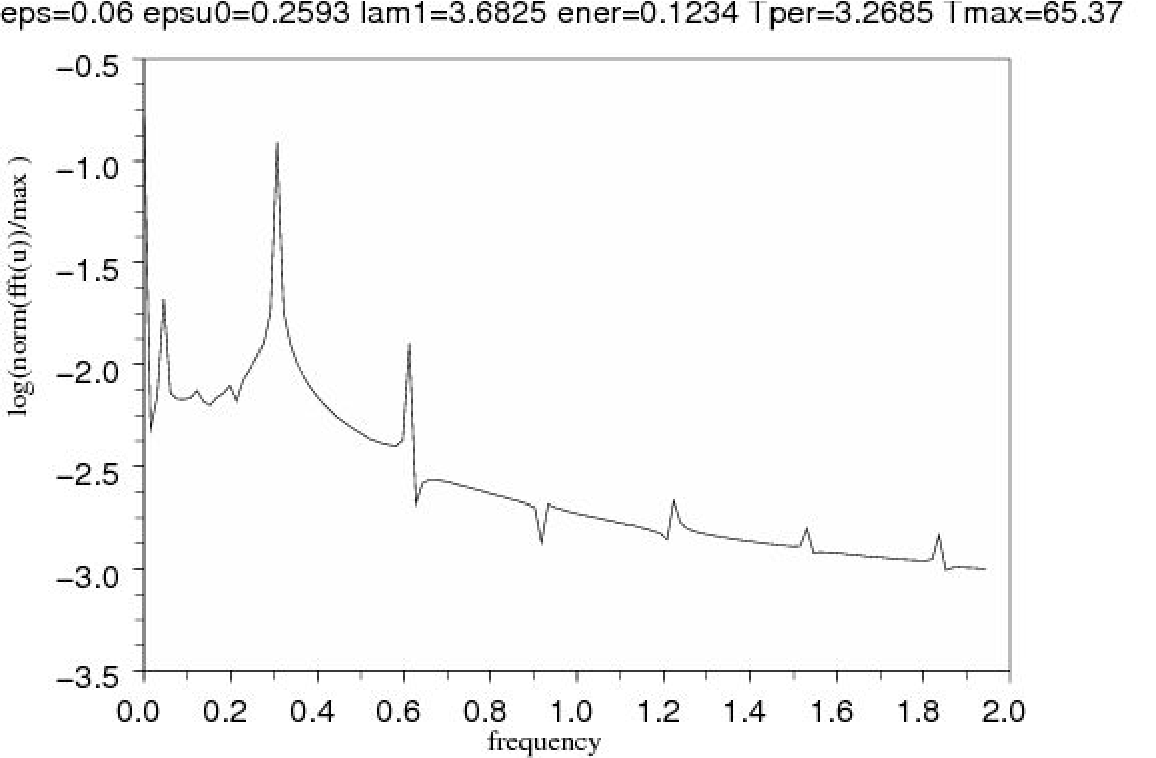}
\includegraphics[ width=7cm, height=7cm]{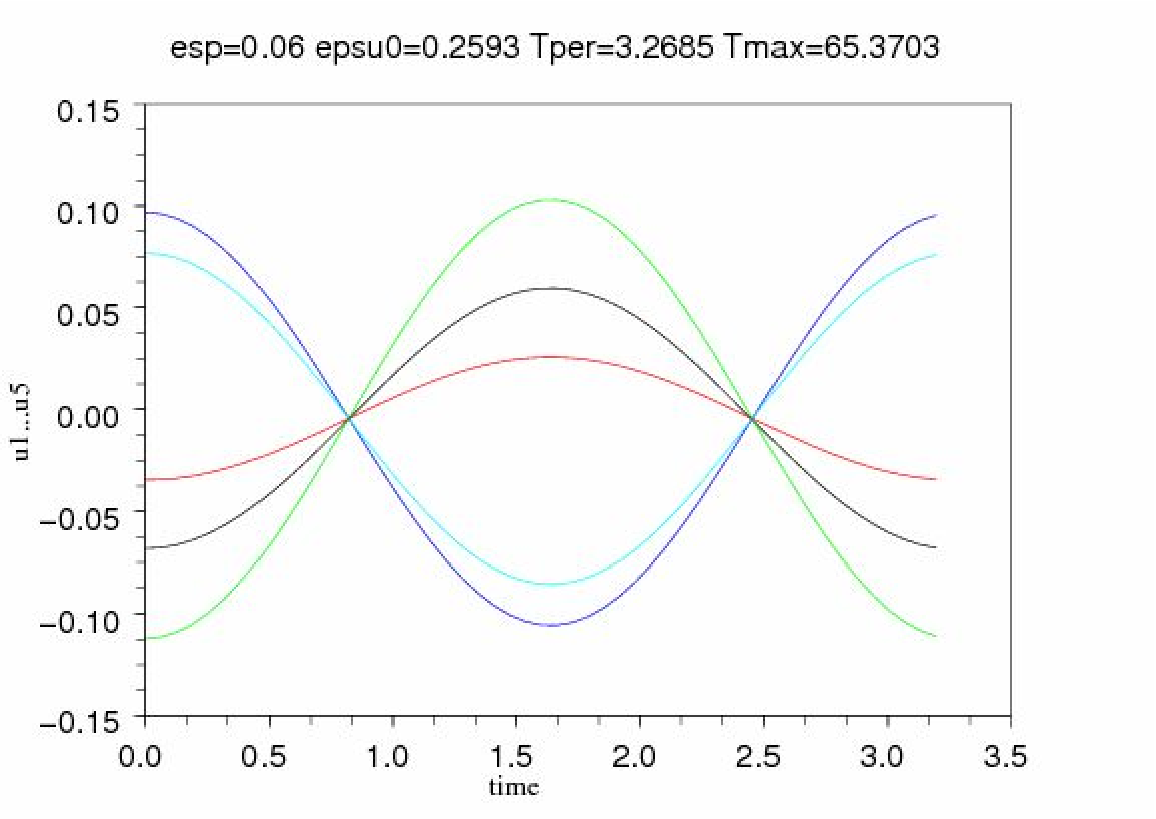}
\caption{energy=0.123, 5 dof; left:Fourier transform ; right:  with respect to time \label{fig:opt-ener0_12}}
\end{center}
\end{figure}

This algorithm may be improved by using  not only the solution associated to the previous value of $\eps$ to solve 
$
 \mathcal{F}(\eps,X_0,X_1,T)=0  
$
but also the derivative of the solution with respect to $X_0,X_1, T$.
\paragraph{Numerical results}
These results are obtained by solving the differential equation with a step by step numerical approximation of the routine $ode$ of Scilab without prescribing the algorithm. As we are looking for a periodic solution, this numerical approximation may be  certainly improved in precision and computing time by using an harmonic balance algorithm.
In figure \ref{fig:opt-ener0_12}, 
 the same example with 5 degrees of freedom and energy equal $0.123$ 
  is  displayed.

 On the left of figure \ref{fig:opt-ener0_12} we find the decimal logarithm of the absolute value of the Fourier transform of the solution; the Fourier transform is computed with the fast Fourier transform with the routine $fft$ of Scilab; we notice the frequency zero due to the non symmetry of the solution and multiples of the  basic frequency; no other frequency appears; on the right the five components are plotted with respect to time; we still notice the non symmetry.

\begin{figure}
  \centering
\begin{center} 
\includegraphics[ width=7cm, height=7cm]{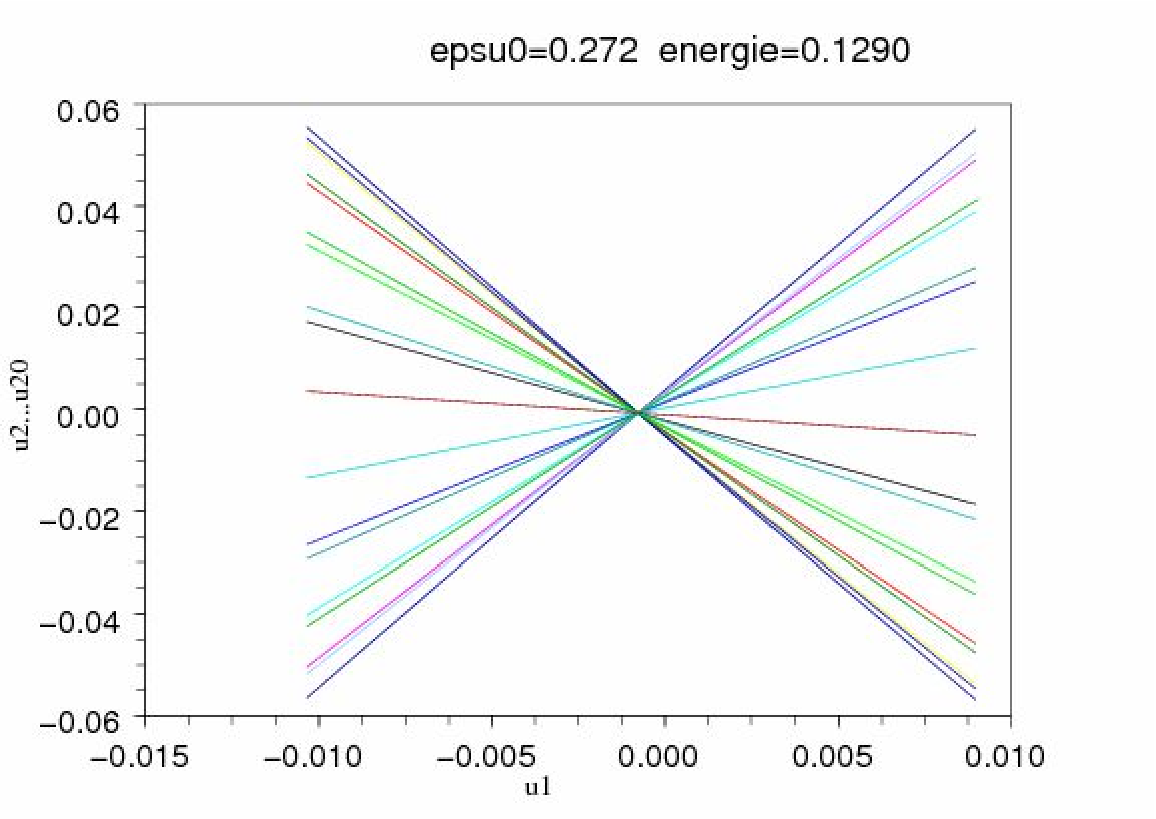} 
  \includegraphics[ width=7cm, height=7cm]{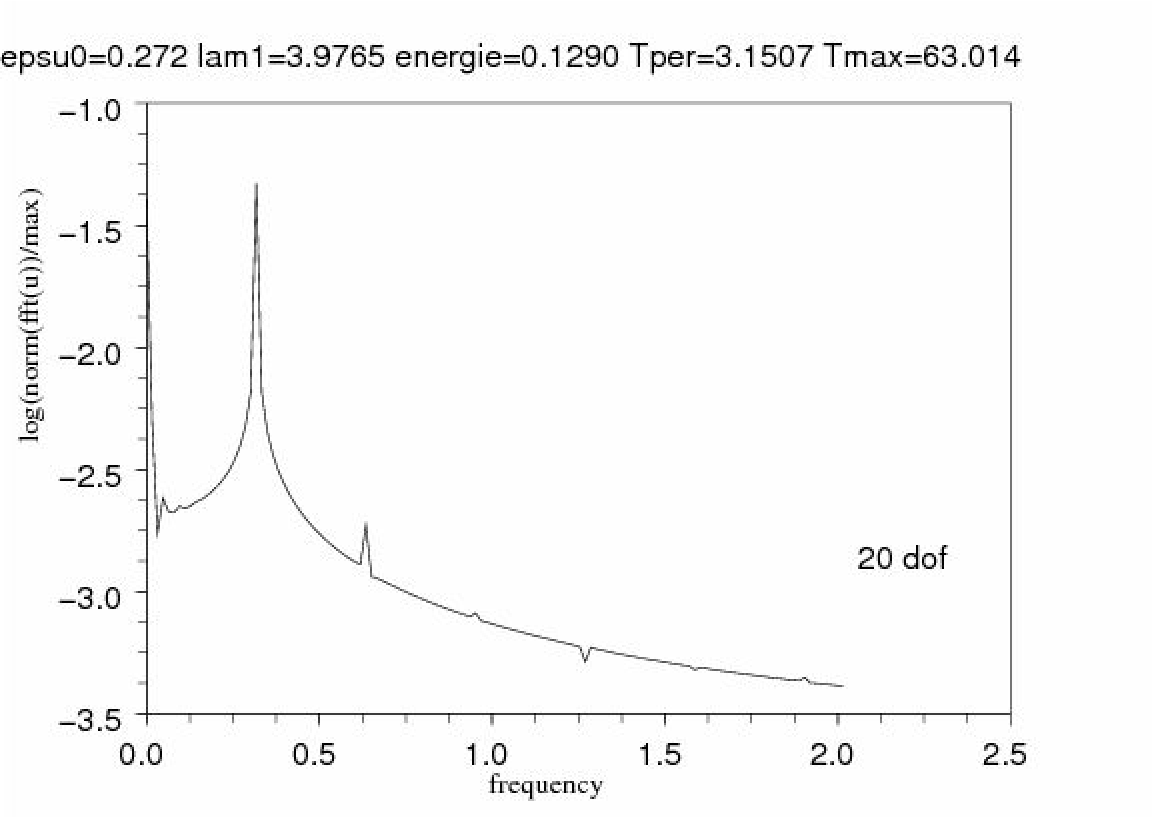}
\caption{energy=0.129, 20 dof; left:in configuration space; right: Fourier transform \label{fig:opt20d-ener0_12}}
\end{center}
\end{figure}

In figure \ref{fig:opt20d-ener0_12} we find results with 20 degrees of freedom, $\eps=0.272$ and energy of $0.129$; the NNM is computed by starting with an  eigenvector associated to the largest eigenvalue . We see on the left in the configuration space that the components are in phase and on the right, the Fourier transform shows zero frequencies and multiple of the basic frequency.

%
\begin{figure}[!ht] 
\centering
\includegraphics[angle=0, width=7cm, height=6cm]{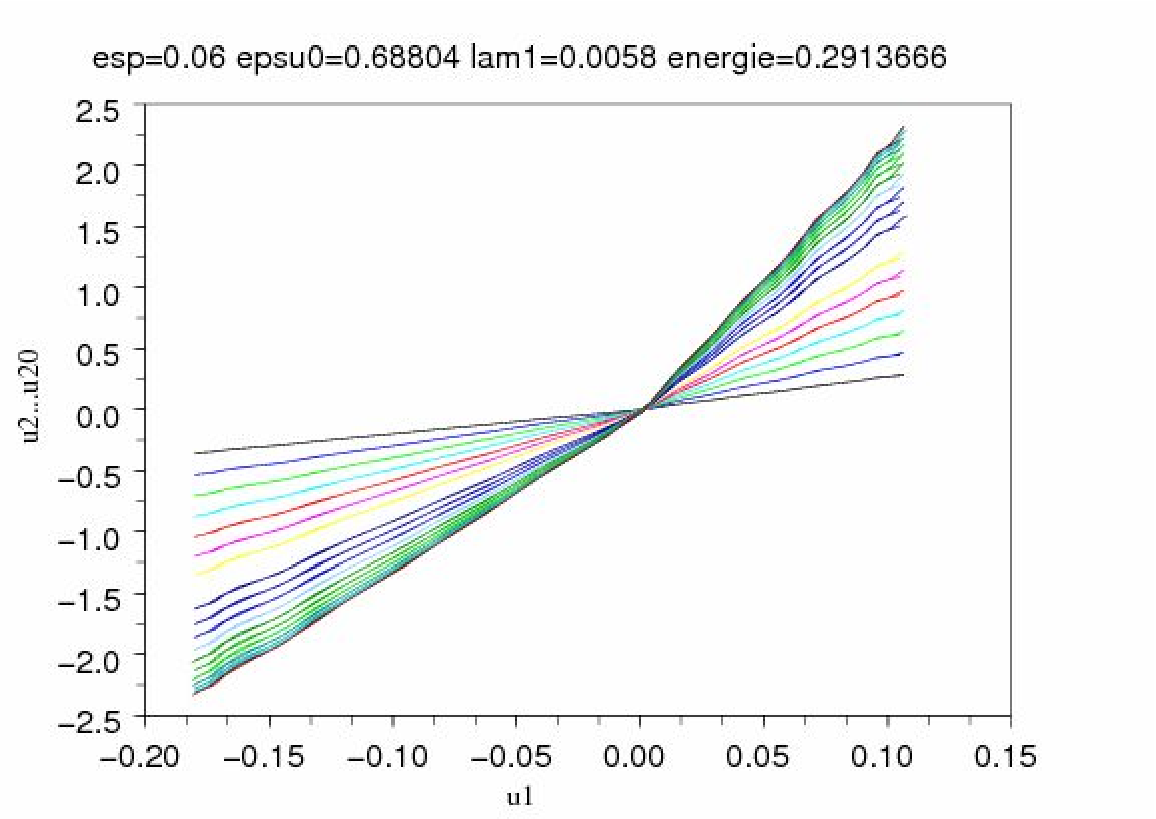}
\includegraphics[angle=0, width=7cm, height=6cm]{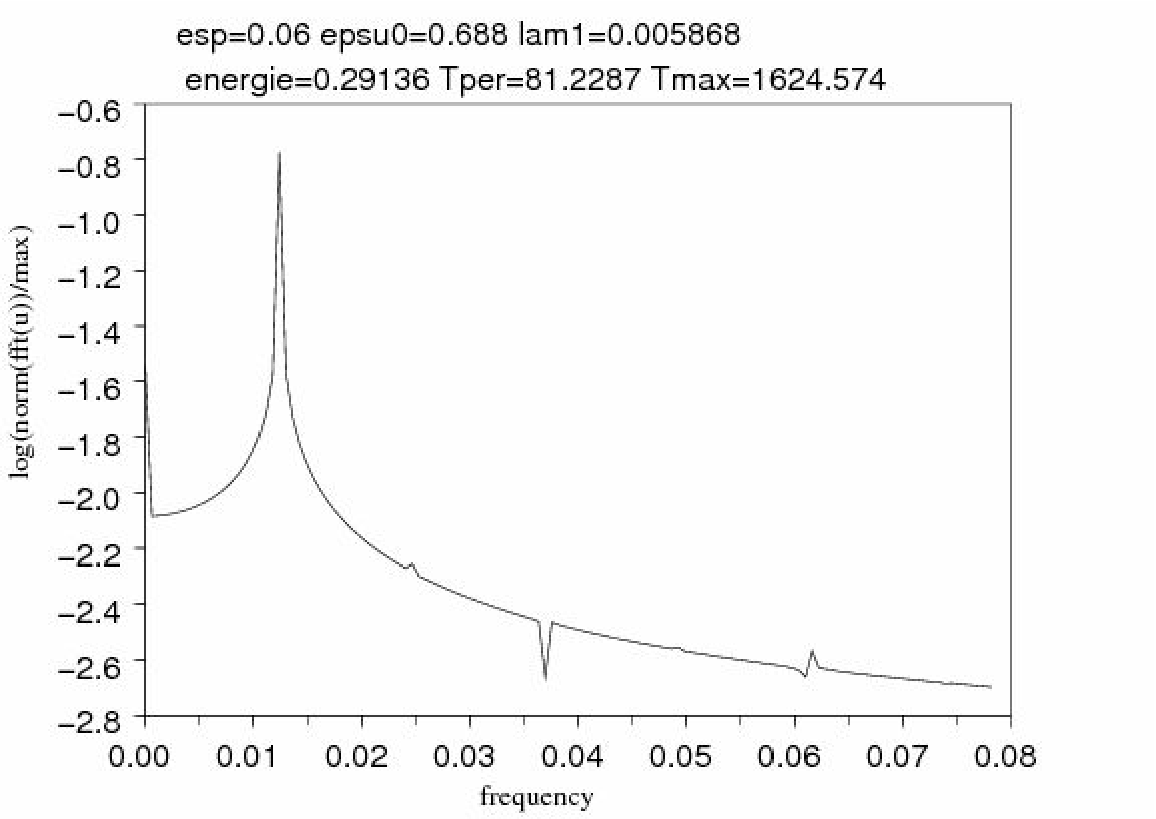}
\caption{ energy=0.29, 5 dof; left: configuration space; right:fft \label{fig:opt20d-ener0_29}}
\end{figure}

In figure \ref{fig:opt20d-ener0_29} the energy is $0.29$ and the NNM is computed by starting with an  eigenvector associated to the smallest eigenvalue; we notice on the left, the solution in the configuration space: at zero each dof has a discontinuity in slope which is clear.

\begin{figure}
  \centering
\begin{center} 
  \includegraphics[ width=7cm, height=7cm]{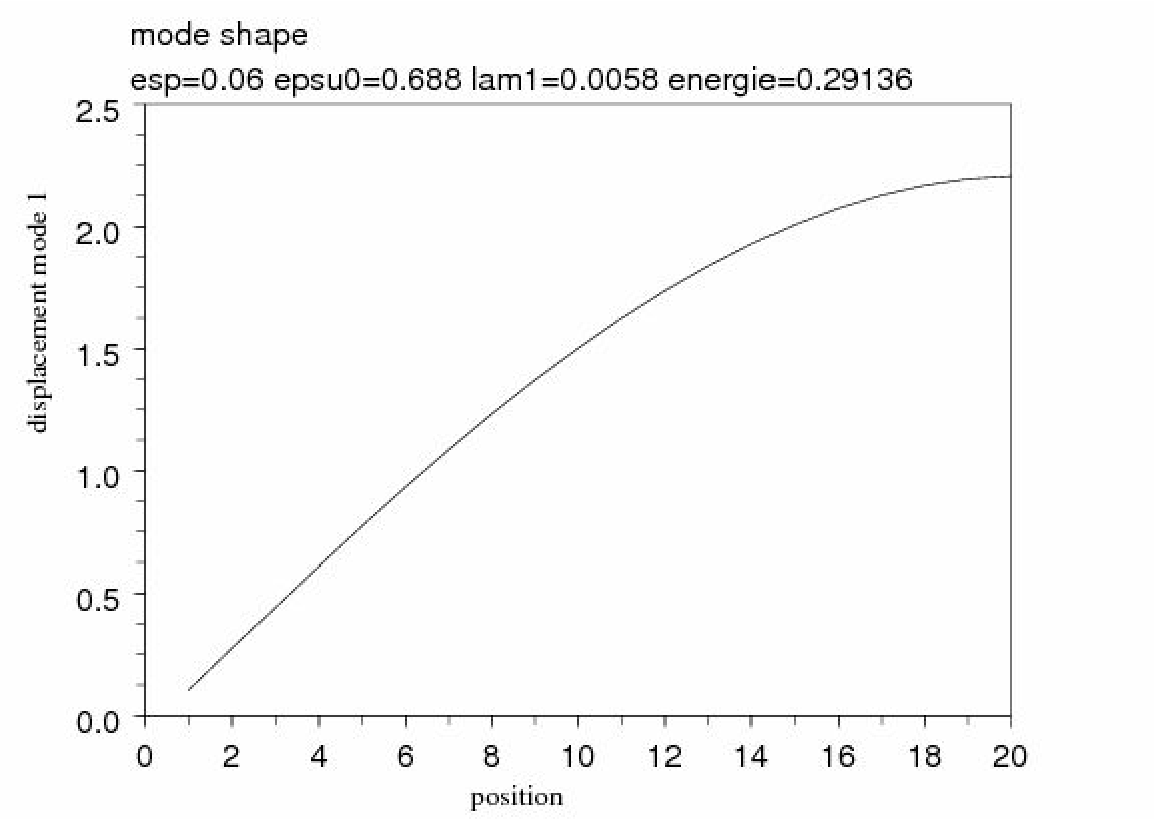}
  \includegraphics[ width=7cm, height=7cm]{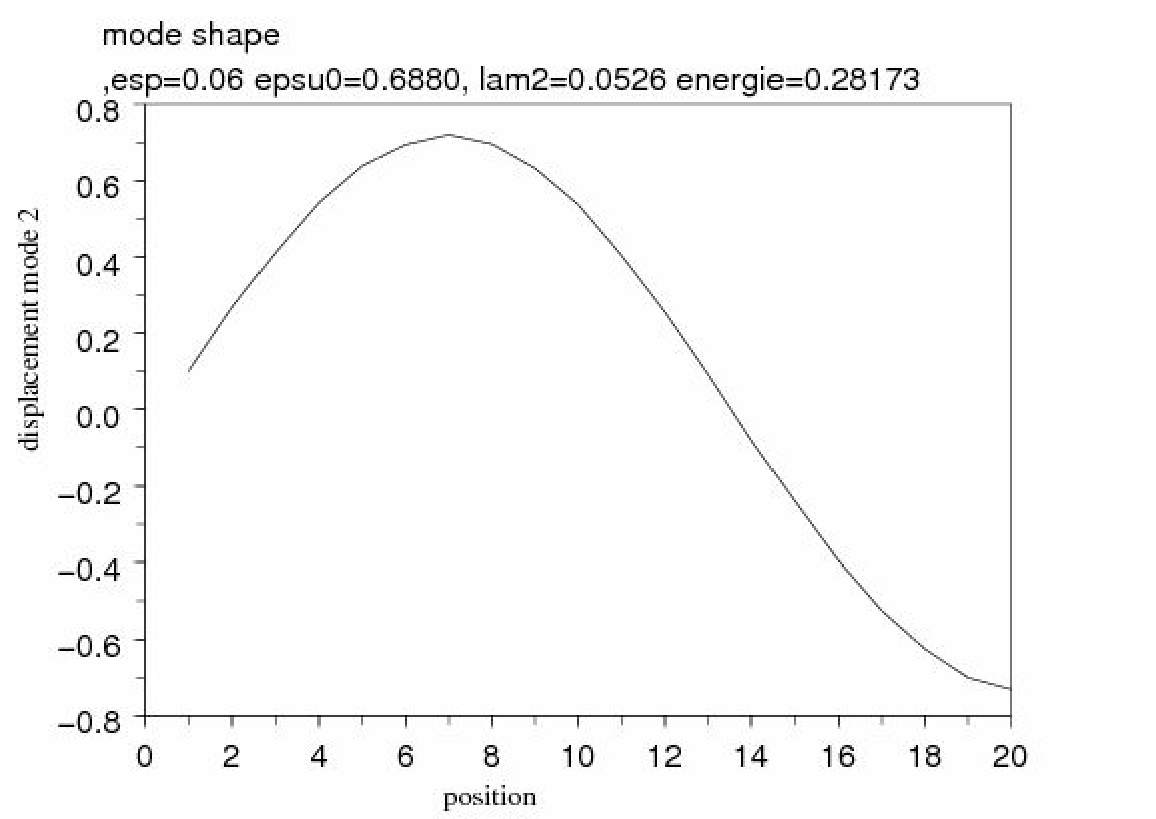}
\caption{ 20 dof; left:energy 0.29 mode 1; right: energy 0.28, mode 2 \label{fig:opt20dmode-ener0_28_29}}
\end{center}
\end{figure}

In figure \ref{fig:opt20dmode-ener0_28_29}, the shape of the NNM is displayed on the left for the NNM  starting from the eigenvector associated to the smallest eigenvalue and on the right for the NNM starting from the second smallest eigenvalue. We notice that the shape is quite similar to the shape of the linear  mode. 




 \subsection{First order asymptotic expansion \label{ssNdof2}}   
%
 In this subsection, we do not particularize the initial data on one eigenmode.
 We adapt  the method  of strained  coordinates when 
  all modes are excited. 
 We loose one order of accuracy compared to previous results  since
 each mode  does not stay periodic and becomes almost-periodic. We assume $\Phi=Id$  to simplify  slightly   the presentation. 
\\
 More precisely,
 the method of strained coordinates is used for each normal component, 
  with the following  initial data
 \bqre  u_k^\eps(0)=a_k,  &\dot{u}_k^\eps(0)=0, & \qquad  k=1,\cdots,N.
  \eqre
 Let us define  $N $ new times $s_k=\lambda_k^\eps t$ and the   following ansatz,
\bqre
  \begin{array}{ccccccl}
   \lambda_k^\eps     &=&  \lambda_k^0 + \eps \lambda_k^1, 
   &&
  \lambda_k^0&=&\lambda_k,
   \\
   u_k^\eps ( t)   & =&   v_k^\eps (\lambda_k^\eps t) &= & v_k^\eps (s_k)
 & =&   v_k^0(s_k) +\eps  r_k^\eps(s_k) . 
 \end{array}
 \eqre
The function $v_k^0$ are easily obtained by the linearized equation.
 Indeed, the only measured nonlinear effect for large time is given 
 by  $(\lambda_k^1)_{k=1}^N$. 
 To obtain these $N$ unknowns, we 
replace  the previous ansatz in the system (\ref{sysMK+}),
 \begin{eqnarray*}
\begin{array}{ccc}
  ( \lambda_k^\eps)^2 (v_k^\eps)''(s_k)  + \lambda_k^2 v_k(s_k)  
  & = & -\eps \dis\left( \sum_{j=1}^N a_{kj}
v_j^\eps\left(\frac{\lambda_j^\eps}{\lambda_k^\eps}s_k \right) -b_k\right)_+.
\end{array}
 \end{eqnarray*}
The right hand side is written in variable $s_k$ instead of 
$\dis s_j$. 
Performing the expansion with respect to  epsilon powers yields
\bqr
 \nonumber
L_k v_k^0&=& ( \lambda_k^0)^2 (v_k^0)''(s_k)  + \lambda_k^2 v_k^0(s_k)  
 = 0 , \\
-L_k  r_k^\eps(s_k)  
 &=&    \dis\left( \sum_{j=1}^N a_{kj}
v_j^0\left(\frac{\lambda_j^0}{\lambda_k^0}s_k \right) -b_k\right)_+
+ 2\lambda_k \lambda_k^1 (v_k^0)''
+R_k^\eps.
\label{eqNdofRemainderH}
 \eqr
Noting that replacing $ v_j^\eps(s_j)$ by 
$v_j^0 \left(\frac{\lambda_j^0}{\lambda_k^0}s_k \right)$  
in (\ref{eqNdofRemainderH})
implies a secular term of the order $\eps t$,
since  
 $s_j = \dis\frac{\lambda_j^0}{\lambda_k^0}s_k  +\mathcal{O}(\eps t)  $, 
 the functions $v_j^0$ are smooth 
and the map $S \rightarrow S_+$ is one-Lipschitz. 
 These new kind of  errors $\mathcal{O}(\eps t)$ are contained   
in   the remainder of each right hand side:
\bqr
  \label{RFose}
    R_k^\eps(t)  & = & \mathcal{O}(\eps t)+ \mathcal{O}(\eps |r^\eps|), 
 \qquad   |r^\eps| = \dis \sqrt{ \sum_{k=1}^N(r^\eps_k)^2}.
\eqr 
If  $b_k=0$,
we identify the secular term 
with the Lemma \ref{LSp} below and the relation   $S_+=\dis S/2+|S|/2 $. 
Then, we
remove  the resonant term
in the source term for the remainder $r_k^\eps$, which gives us 
 $ \lambda_k^1=\dis \frac{a_{kk}}{4\lambda_k}$. 
\\
If $b_k\neq 0$, we compute $\lambda_k^1$ numerically
 with the following orthogonality condition to $\cos(s)$  written in the 
 framework of almost periodic functions,
\bqre
0 & = & 
  \lim_{T \rightarrow \infty} \frac{1}{T}\int_0^T 
\left[
 \dis\left( \sum_{j=1}^N a_{kj}
v_j^0\left(\frac{\lambda_j^0}{\lambda_k^0}s_k \right) -b_k\right)_+
+ 2\lambda_k \lambda_k^1 (v_k^0)''
\right]\cdot \cos(s)ds .
\eqre
The accuracy of the asymptotic expansion depends on the behavior of the 
 solution $\phi=(\phi_1,\cdots,\phi_N)$ of the $N$ following decoupled 
linear equations with  right coefficients $\lambda_k^1$ to avoid resonance
\bqr \label{linsysacc}
-L_k  \phi_k(s_k)  
 &=&    \dis\left( \sum_{j=1}^N a_{kj}
v_j^0\left(\frac{\lambda_j^0}{\lambda_k^0}s_k \right) -b_k\right)_+
+ 2\lambda_k \lambda_k^1 (v_k^0)''.
\eqr
Furthermore  each function $r_k^\eps$ depends  on all times $s_j$, $j=1,\cdots,N$ 
and becomes almost-periodic, 
i.e.  $ r_k^\eps =r_k^\eps( s_1, \cdots,s_N)$.
 Thus the method of strained
 coordinates, 
 only working for periodic functions, 
fails to be continued. 
\\
 Nevertheless, we obtain the following result proved in the Appendix. 
\begin{theorem}[All modes]\label{Th43} $\mbox{ }$ \\
 If $\lambda_1, \cdots , \lambda_N$ are $\zit$ independent,
 then,  for any $T_\eps= o(\eps^{-1})$, i.e.   
   such that 
 \bqre \dis \lim_{\eps \rightarrow 0} T_\eps = +\infty, 
  & \mbox{ and } &
   \dis \lim_{\eps \rightarrow 0} \eps \times T_\eps = 0,
 \eqre
 we have  for all $k=1,\cdots,N$, 
\begin{eqnarray*}\label{aao1}
        \dis \lim_{\eps \rightarrow 0} \|  u_k^\eps(t)-
                               v_k^0\left(  \lambda_k^\eps \;  t\right) 
                               \|_{W^{2,\infty}(0,T_\eps)} = 0
 \end{eqnarray*}
where 
 $ \lambda_k^\eps=  \lambda_k + \eps \lambda_k^1, $
 $v_k^0(s)= a_k \cos(s)$,  and $\lambda_k^1$ is defined by:
 \bqre  \dis    \lambda_k^1 &=&
\frac{1}{2 \lambda_k a_0}
\lim_{T \rightarrow + \infty}\frac{1}{T} \int_{0}^{ T}  
  \left( \sum_{j=1}^N a_{kj}
v_j^0\left(\frac{\lambda_j^0}{\lambda_k^0}s_k \right) -b_k\right)_+
 \cos(s) ds.   
\eqre
Furthermore, if  $b_k=0$, the previous integral yields: 
$\dis \lambda_k^1=\frac{a_{kk}}{4\lambda_k}$.
\end{theorem}
Notice that accuracy  and large time  are weaker  than
 these obtained in Theorem \ref{Th41}. 
It is due to  the inevitable accumulation of the spectrum near the resonance 
and the various times using in the expansion. 
On the other side we have the following direct 
 improvement from  the Theorem \ref{Th41}: 
\begin{remark}[Polarisation] \label{polarisation}
 If only  one mode are excited, for instance the number 1, 
 i.e.  $ a_1 \neq 0$,  $a_k=0$ for all $k\neq1$, then we have the estimate
 for all  $ t \in [0,\eps^{-1}]$:
\bqre 
 \dis   u_1^\eps(t)=&    v_1^0\left(  \lambda_1^\eps \;  t\right) 
                                &+  \mathcal{O}(\eps). 
                               \\
         u_k^\eps(t)=&   0
                &+ \mathcal{O}(\eps) \qquad \text{ for all } k \neq 1.
\eqre  
\end{remark}

    \section{Expansions with even periodic functions
               \label{scos}} 

   Fourier expansion  involving only cosines are used throughout this paper.
 There is never sinus.  
  In this short section we explain  why it is simple to work with even periodic
 functions 
 and we give some hints  to work with  more
 general initial data.  
\\

 First, we want to work only with co-sinus to avoid two secular terms.
 If we return to equation (\ref{eq1onedof}): 
$ -\alpha_0  ( v_1''  + v_1)  =   (v_0)_+ +  \alpha_1 v_0''$.
 A priory, we have two secular terms in the right hand side,
 one with $\cos(s)$ and another
 with $\sin(s)$. Only one  parameter $\alpha_1$ seems not enough to cancel out all secular terms. 
\\
 Otherwise, if $ v_0 \in \rit$, 
 $u, \, S $ are  $2\pi$  periodic even functions,
$g \in C^0(\rit,\rit)$  such that
\\ 
$\dis 
 0=  \int_0^{2\pi} e^{is} (  S(s)+g(u(s))) ds 
$
 then  the solution of
\bqre 
  v''  + v  &= & S(s)+g(u), \quad v(0) =v_0, \, v'(0)=0,
\eqre
 is necessarily a  $2\pi$  periodic even function. 
  Since  we only work with  $2\pi$  periodic even functions
  we have always at most one secular term 
 proportional to  $\cos(s)$.

Now we investigate the case involving not necessarily even periodic functions.
 In general, $\dot{u}_0^\eps \neq 0 $ and
   $u_\eps$ is the solution of 
\bqre
\ddot{u_\eps}  + u_\eps+\eps f(u_\eps)  =  0,
  \quad u_\eps(0) =u_0^\eps, \, \dot{u_\eps}(0)=\dot{u}_0^\eps.
\eqre
By the energy $2E = \dot{u}^2 + u^2 + \eps F(u)$, where $F'=2f$ and $F(0)=0$, 
we know  
 that $u_\eps$ is periodic for $\eps$ small enough, 
 for instance  with an implicit function theorem see \cite{Ver}
 also valid for Lipschitz function \cite{Cla} in our case.  
Denote by $\tau_\eps$
 the first time such that $\dot{u}_\eps(t)=0$. Such time exists 
 thanks to the periodicity of $u_\eps$.  Now, let $U_\eps$ 
 defined by $U_\eps(t)= u_\eps(t+\tau_\eps)$. $U_\eps $ is the solution 
 of  
 \bqre
\ddot{U_\eps}  + U_\eps+\eps f(U_\eps)  =  0,
  \quad U_\eps(0) =U_0^\eps=u_\eps(\tau_\eps), \, 
 \dot{U_\eps}(0)=0.
\eqre
The initial data $U_0^\eps$ depends on the initial position and 
 initial velocity of $u_\eps$ through the energy, 
$(U^\eps_0)^2 + \eps F(U_0^\eps)= (u_0^\eps)^2 + (\dot{u}_0^\eps)^2
   + \eps F(u_0^\eps) $.  For instance, if $u_0^\eps$
 and $\dot{u}_\eps^0$ are positive then $U_0^\eps $ is positive 
 and
\\ 
$U_0^\eps= \dis \sqrt{(u_0^\eps)^2 + (\dot{u}_0^\eps)^2  
  + \eps (F(u_0^\eps)- F\left(\sqrt{(u_0^\eps)^2 + (\dot{u}_0^\eps)^2} \right) }
 + \mathcal{O}(\eps^2)  
          $ .
\\
We can apply the method of strained coordinates for $U_\eps$ 
only with  even periodic functions: 
  $U_\eps (t)= v_0(\omega_\eps t) 
 + \eps  v_1(\omega_\eps t) +\mathcal{O}(\eps^2)$.
The expansion obtained for $u_\eps$ by $U_\eps$, 
with $\phi_\eps =-\omega_\eps \tau_\eps$  is:
 \bqre 
 u_\eps(t)& =&   v_0(\omega_\eps t +\phi_\eps ) 
 + \eps  v_1(\omega_\eps t +\phi_\eps) +\mathcal{O}(\eps^2),
\eqre
   which is a  good ansatz in general for $u_\eps$, 
  where $v_0$ and $v_1$ are even $2\pi-$periodic functions. 
 The method of strained coordinates becomes to find the following 
 unknowns  $\phi_0$, $\omega_1, \phi_1$, $\omega_2, \phi_2$ such that 
 \bqre
  \omega_\eps &= &\omega_0 +\eps \omega_1 +\eps^2 \omega_2 + \cdots,\\
   \phi_\eps &=& \phi_0 +\eps \phi_1 +\eps^2 \phi_2 + \cdots.
\eqre
 Indeed,  we have two parameters to cancel out  two secular terms at each step.
 If one is  only interested by the nonlinear frequency shift, it is simpler 
 to work only with cosines. 
 \\
   
 Otherwise, if $f$ is an odd function, we can work only with 
 odd periodic function. It is often the case in literature when occurs a cubic
 non-linearity. See for instance \cite{KC96,Mi,Nayfeh} for the Duffing equation,
 the Rayleigh equation or  the Korteweg-de Vries equation.

    \section{Appendix: technical proofs \label{sA}} 

We  give some useful results  about energy estimates and almost periodic functions in subsection
\ref{ssul}.
 Next we  complete the  proofs for each previous asymptotic expansions
 in subsection \ref{ssbr}. 
 The   point is to bound the remainder for large time
 in each expansion.
\subsection{Useful lemmas} \label{ssul}
The following Lemma is useful to prove an expansion 
for large   time with non smooth non-linearity.
%
\begin{lemma}{\bf [Bounds for large time 
   ]}
\label{lemmaTime}
 $\mbox{ }$
\\
 Let   $w_\eps$ be a    solution of  
 \bqr \label{eqGeneralRemainder} 
 \dis \left\{
\begin{array}{l}
  w_\eps''+  w_\eps   
 =  S(s) +  f_\eps(s) + \eps g_\eps(s,w_\eps),
    \\
   w_\eps(0)=0, \quad  w_\eps'(0)=0.
 \end{array} 
\right.
 \eqr
  If  source terms satisfy the following conditions where  
 $M > 0$, $C >0$ are fixed constants :
\begin{enumerate}
\item 
 $S(s)$ is a $2\pi$-periodic function orthogonal to 
  $e^{\pm i s}$,
and   $|S(s)|\leq M$ for all $s$,
  \item 
   $|f_\eps|\leq M $ and for all $T$,
  $\dis \int_0^T |f_\eps(s) | ds  \leq C\eps T \quad (\mbox{resp. } C\sqrt{\eps}T$),
 \item
 for all $R> 0$: 
 $ \dis M_R=\sup_{\eps \in (0,1),s >0, R > |u|} |g_\eps(s,u)| < \infty$, 
\\ that is to say that  
  $g_\eps(s,u)$ is locally bounded  with respect to $u$ 
  \\ 
 for $\eps \in (0,1)$ and $s \in (0,+\infty)$,
\end{enumerate}
then,   there exists $\eps_0 > 0$  and $\gamma > 0$
    such that, for $ 0 < \eps < \eps_0$,  
   $ w_\eps $ is uniformly bounded in $ 
         \dis W^{2,\infty}\left(0,T_\eps \right)$,
  where $T_\eps = \dis \frac{\gamma}{\eps} 
  \quad (\mbox{resp. } \dis \frac{\gamma}{\sqrt{\eps}}) $.
\end{lemma} 
 Notice that $f_\eps$ and $g_\eps$ are not necessarily continuous.
 Indeed this a case for our asymptotic expansion, 
 see Lemme \ref{LDLL} and its applications throughout the paper.
But in previous sections the right hand side is globally continuous, 
 i.e. $S+ f_\eps+\eps g_\eps(.,w_\eps)$ is continuous, so, in this case, 
 $w_\eps$ is $C^2$. 
\\
\medskip
\\
{\bf Proof of the  Lemma \ref{lemmaTime}}: 
 First we remove the non resonant periodic source term which is independent of $\eps$. Second,
 we get $L^\infty $ bound for $w_\eps$ and $w_\eps'$ with an energy estimate.
 Third, with equation (\ref{eqGeneralRemainder}), we get 
 an uniform estimate for 
 $w''_\eps$ in $L^\infty(0,T_\eps)$ and the $W^{2,\infty}$ regularity.
\\
\underline{Step 1: remove $S$}
 \\
  It suffices to  write $w_\eps = w_1 + w_2^\eps$ where 
  $ w_1$ solves the linear problem: 
  \bqr 
  \label{eqw1}  
 w_1''+  w_1   =  S(s),
    \quad   w_1(0)=0,  w_1'(0)=0. 
 \eqr
$w_1$ and $w_1'$ are uniformly bounded 
 in   $L^\infty(0,+\infty)$
  since there is no resonance.
\\
More precisely, 
 $w_1 = F(s) 
 + A \cos(s) + B \cos(s)$,
 where $F$ is  $2\pi$ periodic.
$ F$ is obtained by Fourier expansion without harmonic $n=\pm 1$
 since $S$ is not resonant:
\bqre 
 F(s) &= &
 \sum_{ n \neq \pm 1} 
  \frac{c_n}{1 - n^2}e^{ins}
\quad 
 with
  \quad 
  S(s) = 
 \sum_{ n \neq \pm 1} c_n e^{ins}.
\eqre
$F$ is uniformly bounded, 
 with Cauchy-Schwartz inequality set 
 $C_0^2 = \dis \sum_{ n \neq \pm 1} 
  |n^2-1  |^{-2}$, 
 we obtain: 
$
 \|F \|_{L^\infty}  \leq  
 \sum_{ n \neq \pm 1} 
  \frac{|c_n|}{|n^2-1|}
\leq  C_0  \|S \|_{L^2(0,2\pi)}
\leq  C_0  \|S \|_{L^\infty(0,2\pi)}.
$
\\
Similarly, set   $D_0^2 =  \sum_{ n \neq \pm 1} 
 n^2 |n^2-1 |^{-2}$, 
we have 
 $  \|F' \|_{L^\infty} \leq  D_0  \|S \|_{L^\infty(0,2\pi)}$.
\\
Furthermore, $0=w_1(0)= F(0)+ A$, 
  and   $0=(w_1)'(0)= F'(0)+ B$,
 then,  $A$ and $B$  are well defined. 
 $w_1^,$  is also bounded, i.e. 
 there exists $M_1> 0 $ such that  
 $\|w_1\|_{W^{1,\infty}(0,+\infty)} \leq M_1$.
\\
Notice that from equation \eqref{eqw1}, $w_1$ belongs to $W^{2,\infty}$. 
\\
\medskip
\\
Then we  get  an equation similar to (\ref{eqGeneralRemainder}) for $w_2^\eps$ 
 with   $S \equiv 0$ and the same assumption for 
  the same $f_\eps$ and the new $g_\eps$:
 $
    \overline{g}_\eps(s,w) =
     g_\eps(s,w_1 +w)  .$
 \bqr \label{eqGeneralRemainder2} 
 \dis \left\{
\begin{array}{l}
  (w^\eps_2)''+ (w^\eps_2)   
 =    f_\eps(s) + \eps \overline{g}_\eps(s,w_2^\eps),
    \\
   (w_2^\eps)(0)=0, \quad  (w_2^\eps)'(0)=0.
 \end{array} 
\right.
 \eqr
\\
\underline{Step 2: energy estimate}
 \\
 Second, we get an energy estimate for $w_2^\eps$.
   We fix $R> 0$ such that $ R $ is greater than the uniform bound
 $M_1$ 
 obtained for $w_\eps^1$ 
 and yet  $ R= M_1 + \rho$ with $\rho > 0$.
Let us define  
\bqre
  2E(s)=((w_2^\eps)'(s))^2 +  (w_2^\eps)(s)^2 ,
 & \quad 
 \dis \overline{E}(s) =\sup_{0<\tau < s} E(\tau) ,
\eqre
 and $T_\eps$ be the first time   $T>0$ such that 
 $ 2\overline{E}(T)\geq \rho^2$, i.e. 
  $\rho$ estimates the size of $(w_2^\eps)$ and $(w_2^\eps)'$. 
\\
Multiplying the differential equation (\ref{eqGeneralRemainder2})
 by $(w_2^\eps)'$, we have  for  all $s < T < T_\eps(\rho)$
 the following inequalities 
since 
 $\dis \sup_{0<\tau < s} |(w_2^\eps)'(\tau)| \leq \sqrt{2\overline{E}(s)}$,
and 
$\dis \int_0^T |f_\eps(s) | ds  \leq C\eps T
$,
\bqre
E(s)  & = &  \dis \int_0^s f_\eps(\tau ) (w_2^\eps)'(\tau) d\tau
                  + 
 \eps \int_0^s \overline{g}_\eps(\tau,(w_2^\eps)(\tau))(w_2^\eps)'(\tau) d\tau, 
         \\    
      &  \leq &
  C \eps s  \sqrt{2\overline{E}(s)}
          +  \eps s  M_R  \sqrt{2\overline{E}(s)}
        ,\\
\overline{E}(T)   &\leq &  C \eps T  \sqrt{2\overline{E}(T)}
          +  \eps T  M_R  \sqrt{2\overline{E}(T)}, \\
\eps T & \geq & \dis  \frac{\sqrt{\overline{E}(T)/2}}{M_R + C}.
\eqre
Notice that if $ 2 \overline{E}(T) < \rho^2 $ for all  $T > 0$
 then $T_\eps= + \infty $. The critical case is  
when $T_\eps$ is finite and 
  ${\overline{E}(T)}$ approaches $\rho^2/2$ when $T$ goes
  to $T_\eps(\rho)$. 
Thus we have 
 $ \dis 
   \dis T_\eps \geq \frac{\rho}{2\eps(M_R+C)}
$ and $E(t) \le \frac{\rho^2}{2}$ for $t \le T_{\eps}=\frac{\gamma}{\eps}$ with
$\gamma=\frac{\rho}{2(M_R+C)}$.
\\
The proof is similar when 
$\dis \int_0^s |f_\eps(\tau )|  d\tau \leq C \sqrt{\eps} T $
 then 
 $ \dis 
   \dis T_\eps \geq \frac{\rho}{2\sqrt{\eps}(\sqrt{\eps}M_R+C)}
$.
\epro 
 For completeness, we state  a similar and straightforward 
version of Lemma \ref{lemmaTime}
 useful for systems. 
\begin{lemma}{\bf [Bounds for large time for systems
                ]}
\label{lemmaTimeSystem}
 $\mbox{ }$
\\
 Let   $w_\eps= (w_1^\eps,\cdots, w_N^\eps)$
 be the     solution of  the following system:
 \bqr \label{eqGeneralRemainderS} 
 \dis \left\{
\begin{array}{l}
 (\lambda_1)^2 (w_k^\eps)''+ (\lambda_k)^2 w_k^\eps   
 =  S_k(s) +  f_k^\eps(s) + \eps g_k^\eps(s;w_\eps),
    \\
   w_k^\eps(0)=0, \quad  (w_k^\eps)'(0)=0,
 \quad k=1,\cdots,N.
 \end{array} 
\right.
 \eqr
  If  source terms satisfy the following conditions where  
 $M > 0$, $C >0$ are fixed constants :
\begin{enumerate}
\item 
      \underline{non resonance conditions}
       with  $S_k(s)$ are $2\pi$-periodic functions and   $|S_k(s)|\leq M$,
\begin{enumerate}
  \item 
    $S_1(s)$ is orthogonal to 
  $e^{\pm i s}$, i.e. $\dis \int_0^{2\pi} S_1(s)e^{\pm i s}ds =0,$
   \item
    $\{ \lambda_2, \cdots, \lambda_N \} \nin \lambda_1\zit$,
 \end{enumerate}
  \item 
   $|f_k^\eps|\leq M $ and for all $T$,
  $\dis \int_0^T |f_\eps(s) | ds  \leq C\eps T \mbox{ or } C\sqrt{\eps}T$,
 \item
 for all $R> 0$: 
 $ \dis M_R=\max_k\sup_{\eps \in (0,1),s >0, w_1^2+\cdots+w_N^2< R^2} 
 |g_k^\eps(s;u)| < \infty$, 
\end{enumerate}
then,   there exists $\eps_0 > 0$  and $\gamma > 0$
    such that, for $ 0 < \eps < \eps_0$,  
   $ w_\eps $ is uniformly bounded in $ 
         \dis W^{2,\infty}\left(0,T_\eps \right)$,
  where $T_\eps = \dis \frac{\gamma}{\eps} 
  \mbox{ or } \dis \frac{\gamma}{\sqrt{\eps}} $.
\end{lemma} 
\bpro
 First we remove source terms $S_k$ independent of $\eps$
setting $w_k^\eps = w_{k,1}+w_{k,2}^\eps$ where $w_{k,1}$ is the solution of 
\bqre
     \lambda_1^2 w_{k,1}'' + \lambda_k^2 w_{k,1}=S_k, & 
    w_{k,1}(0)= 0, &  w_{k,1}'(0)=0.
\eqre
As in the proof of  Lemma \ref{lemmaTime}, 
$w_{1,1}$ belongs in $W^{2,\infty}$ thanks to 
 the non-resonance condition {\it 1.(a)}. 
For $k \neq 1$, there is no resonance since
  $\frac{\lambda_k}{\lambda_1} \nin \zit,$
 i.e.  the non-resonance condition {\it 1.(b)}, 
 thus a similar expansion also yields
$w_{k,1}$ belongs in $W^{2,\infty}(\rit,\rit)$.
\medskip \\ 
Now  $w_{k,2}^\eps$ are solutions of the following system for $k=1,\cdots,N$
 \bqre  
 \dis \left\{
\begin{array}{l}
 \lambda_1^2 (w_{k,2}^\eps)''+ \lambda_k^2(w_{k,2}^\eps)   
 =    f_k^\eps(s) + \eps \overline{g}_k^\eps(s;w_2^\eps),
    \\
   (w_{k,2}^\eps)(0)=0, \quad  (w_{k,2}^\eps)'(0)=0,
 \end{array} 
\right.
 \eqre
with $w_\eps=w_1+w_2^\eps$, $w_2^\eps=(\cdots,w_{k,2}^\eps,\cdots)$ and 
 $
    \overline{g}_k^\eps(s;\cdots,w_k,\cdots) =
     g_k^\eps(s;\cdots,w_{k,1} +w_k,\cdots)  .$
\\
The end of the proof of  Lemma \ref{lemmaTimeSystem}
 is a straightforward generalization of the  
the proof of  Lemma \ref{lemmaTime} with the 
energy:
 $ 2 E(w_1,\cdots,w_N) = \dis 
\sum_{k=1}^N \left( 
 ( \lambda_1)^2 (\dot{w}_k)^2 + (\lambda_k)^2 w_k^2
 \right)
.$
\epro
\medskip
 For systems,  we also have to work 
 with linear combination of periodic functions 
with different periods and nonlinear function of such sum.
 So we  work with the adherence in $L^\infty(\rit,\cit)$ of 
 span$\{e^{i\lambda t},\, \lambda \in \rit\}$, namely 
 the set  of almost periodic functions 
$ C^0_{ap}(\rit,\cit)$, 
and the  Hilbert space of almost-periodic function 
 is  $ L^2_{ap}(\rit,\cit)$,
  see \cite{Cor},
with the scalar product 
\bqre 
 \left <u,v \right > 
 & = &\lim_{T \rightarrow + \infty}
                \frac{1}{T}\int_0^T u(t) \overline{v(t)}dt.
\eqre
 We  give an useful Lemma about the spectrum of $|u|$
     for  $u \in C^0_{ap}(\rit,\rit)$.   
   Let us recall  definitions for the Fourier coefficients 
  of $u$ associated to frequency $\lambda$: $ c_{\lambda}[u]$
 and its spectrum: $Sp\;[u]$,
 \bqr 
 \label{psap}
  \dis
 c_{\lambda}[u] =  \left < u,e^{i\lambda t} \right > 
                  = \lim_{T \rightarrow + \infty}
                \frac{1}{T}\int_0^T u(t) e^{-i\lambda t}dt,
 & \quad &
 Sp\;[u] = \{\lambda \in \rit, \, c_\lambda[u] \neq 0\}.
\eqr
 \begin{lemma} \label{LSp}
 {\bf [Property of the  spectrum of $|u|$  ] }  \alali
Let   $u \in C^0_{ap}(\rit,\rit)$ a function with  a finite spectrum: 
$Sp\;[u] \subset \{\pm \lambda_1,\cdots, \pm \lambda_N\}$.
\\
If
$ (\lambda_1,\cdots,\lambda_N)$ are  $\zit$-independent,
then 
  $\lambda_k \nin Sp\;[\; |u| \;]$  for all $k$. 
\end{lemma}
\bpro 
  Notice that  $ 0 \nin Sp\;[u]$.
   The result is quite obvious for $u^2$.  
We first prove the result for $f(u^2)$ where $f$ is smooth.
 Then, we conclude by approximating $|u|$ by a smooth sequence 
$f_n(u^2)=\sqrt{1/n+u^2}$, 
  and using the $L^\infty$ stability of the spectrum.
\\
Let $E$ be the set of all $\zit$ linear combinations of elements of
$ S_2 = \{0, \pm \lambda_{kj}^\pm, \,  k,j=1,\cdots N\}$, 
where 
 $ 
   \lambda_{kj}^\pm = \lambda_k\pm \lambda_{j}.
$
 Thus $ Sp\;[f(u^2)]$ is a subset of $E$   since   $Sp\;[u^2]  \subset S_2$.
\\
Notice that 
$ 
 \lambda_{jk}^\pm=\pm  \lambda_{kj}^\pm,
           \quad
             \lambda_{kk}^{-} =0,
 \quad
  \lambda_{kk}^+ =2\lambda_k=\lambda_{kj}^++ \lambda_{kj}^-.
$ \\
 Choosing $k=1$ for instance, so $\lambda_1 \neq 0$, it suffices to prove
 that $\lambda_1 \nin E$.   
\\
Assume the converse, 
   i.e., $\lambda_1 \in E$. 
 Then, for $k<j$,
  there exists some integers 
 $\dis (c_{kj}^\pm)_{k < j)}$ such that: 
$ 
 \lambda_1  =  \sum_{k < j}( c_{kj}^+ \lambda_{kj}^+ +  
                         c_{kj}^- \lambda_{kj}^-).
$ 
 Therefore, defining $c_{jk}^\pm$ by  $\pm c_{kj}^\pm$ for $k < j$, 
  we have:
\bqre 
   \lambda_1  &= & \lambda_1  \sum_{ j\neq 1}( c_{1j}^+  + 
                         c_{1j}^-) 
         + \lambda_2  \sum_{ j\neq 2}( c_{2j}^+  
                        +  c_{2j}^-) 
       + \cdots + \lambda_N\sum_{ j\neq 1}( c_{k=Nj}^+  
                         + c_{Nj}^-). 
\eqre
Using the $\zit$-independence, with $ d_{kj}= c_{kj}^+  + 
                         c_{kj}^-$ for $k \neq j$ and $d_{kk}=0$, we have 
 following system:
$\dis 
  1  =  D_1= \sum_j d_{1j}, \quad 
  0  =  D_k= \sum_j d_{kj}, \qquad  \mbox{for all  } k >1. 
$
\\
Summing up, the $N-1$ last equations  in $\dis \frac{\zit}{2\zit}$,
 and using the fact: $ d_{jk} \equiv d_{kj}$ modulo  $2$,
 we have:
$
\dis 
  0 \equiv \sum_{k=2}^{N} D_k
   \equiv \sum_{j=2}^{N} d_{1j} + 2\sum_{k<j}^{} d_{kj} 
      \equiv  \sum_{j=2}^{N} d_{1j},
$
then 
$ \dis D_1  \equiv 0$, i.e.  $D_1$ is even.
 It's impossible since $D_1=1$. 
   So $\lambda_1 \nin E $ 
and the proof is complete.
 \epro

\subsection{Bounds for the remainders}\label{ssbr}

Now, we prove each asymptotic expansion given in previous sections, i.e. 
 we bound each remainders with energy estimates up to a large time.
\medskip 
\\
{\bf Proof of  Proposition \ref{propexplicit} :}
 First we give the outline of the proof.
\\
 Notice that all these computations only involve 
 the function $\cos$. Then, the only way to have a secular term 
in equations defining  $v_1$ and $v_2$ is  a $\cos(s)$
 in the right-hand side. So, the good choice of $\alpha_1$ and $\alpha_2$,  
 is enough to remove secular term with $\cos(s)$. 
Now, it suffices to control $r_\eps$ for large time.
A computation  shows that
the remainder  $R_\eps$ of equation (\ref{rode}) satisfies:
\\ 
$|R_\eps(s)|\leq  C \eps( 1 + |r_\eps(s)|)  
                 +| \chi_\eps|(v_0,v_1+\eps r_\eps)$.  
\\
Then, $r_\eps$ is like $w_\eps$ 
 in Lemma \ref{lemmaTime}, and the term $f_\eps$ comes from 
 $\chi_\eps$ which  is estimated by
%
  Lemmas \ref{LDLL}, \ref{LHzi}.
\\
More precisely, an exact computation of $R_\eps$ in equation (\ref{rode})
leads to 
\bqre
   R_\eps &= &
  \chi_\eps(v_0,v_1+\eps r_\eps) + \eps H(v_0) r_\eps 
 + \eps \alpha_3^\eps v''_\eps, 
\eqre
where 
  $\alpha_3^\eps$ is  a real constant, bounded uniformly for all 
$\eps \in [0,1]$  such that
\\
 $
 (\omega_\eps)^2
 = 
 \alpha_0 + \eps \alpha_1+\eps^2 \alpha_2+\eps^3 \alpha_3^\eps. 
$ 
From (\ref{eqchi}) we also have
\bqre
   & \chi_\eps(v_0,v_1+\eps r_\eps) \\
  = &
 \dis  \left\{ (v_0 +\eps v_1+\eps^2 r_\eps )_+
           - [(v_0)_+ +  \eps H(v_0)( v_1+\eps r_\eps)  ] \right\}\eps^{-1} 
 \\
   =& 
  \left\{(v_0 +\eps v_1 )_+
           - [(v_0)_+ +\eps  H(v_0) v_1  ] 
  +
 (v_0 +\eps v_1+\eps^2 r_\eps )_+
           - (v_0 +\eps v_1)_+   \right\}\eps^{-1} 
 - \eps H(v_0) r_\eps
 \\
=&
 \chi_\eps(v_0,v_1) -\eps H(v_0) r_\eps + \eps \tilde{g}_\eps(s,r_\eps),
\eqre  
 since $ u \rightarrow (u)_+$ is 1-Lipschitz
 $ |\tilde{g}_\eps(s,r_\eps) |  = \left\{(v_0 +\eps v_1+\eps^2 r_\eps )_+
           - (v_0 +\eps v_1)_+   \right\}\eps^{-2}  \leq 
 |r_\eps|.
$
So, with  $v_\eps = v_0+\eps v_1+\eps^2 r_\eps$,
 we can rewrite $R_\eps$ as follow
\bqre 
  R_\eps &= &
  \chi_\eps(v_0,v_1) + \eps  \tilde{g}_\eps(s,r_\eps) 
 + \eps \alpha_3^\eps v''_\eps. 
\eqre
Now, we can rewrite equation (\ref{rode}) in the following way
\bqre 
- \alpha_0(r''_\eps +  r_\eps)
 & =&
 S(s) +f_\eps(s) + \eps g_\eps(s,r_\eps), 
\eqre
with  $ S  =  \alpha_2 v_0'' + \alpha_1 v_1''  + H(v_0) v_1,
      $
 $ f_\eps  =  \chi_\eps(v_0,v_1)
                  + \eps\alpha_3^\eps(v_0"+\eps v_1"), 
 $
$
 g_\eps = \tilde{g}_\eps 
                   + \eps^2\alpha_3^\eps r_\eps,
$
which allows us to conclude with Lemma \ref{lemmaTime}.
\epro 
The proof for other propositions \ref{PropNdaf}, \ref{Propcritic} 
 in section \ref{s1dof} are similar. 
\\
We now complete the proof for the asymptotic expansions for systems
given in section \ref{sNdof}.
\medskip \\
{\bf Proof of Theorem \ref{Th41} :}
 As in the proof of  Proposition \ref{propexplicit}, 
 the same technique is used component by component for 
 Theorems \ref{Th41}, 
 with   similar  energy estimates  
 we  can conclude 
 with the  Lemma \ref{lemmaTimeSystem} for system to control all 
 $r_k^\eps$.
\\
To  simplify the writing of the proof, let us  assume that $\Phi=Id$
in \eqref{sysMK+}. In this case, a complete computation of the remainder gives us:
\bqre 
 R_k^\eps &= &
  \chi_\eps(a_{k1}v_1^0 - b_k,\sum_j a_{kj}
[v_{kj}^1+\eps r_j^\eps] ) + \eps H(a_{k1}v_1^0 - b_k)\sum_j a_{kj} r_j^\eps 
 + \eps \alpha_3^\eps (v_k^\eps)'', 
\eqre
 with notation of the  proof of Proposition \ref{propexplicit} and 
 $v_k^\eps = v_k^0 +\eps v_k^1+\eps^2 r_k^\eps$.
Let $u,v,w$ be three functions, as previously, we have:
\bqre
\chi_\eps( u,v+\eps w)+ \eps H(u)w &= &
 \chi_\eps( u,v) + \eps^{-1}( (u+\eps v +\eps^2 w)_+ -(u+\eps v )_+ ) , 
\eqre
 and, since $ w \rightarrow w_+$ is 1-Lipschitz:
$ 
|\eps^{-1}( (u+\eps v +\eps^2 w)_+ -(u+\eps v )_+ )|
 \leq   \eps |w|.
$
Now, we can rewrite $R_k^\eps$ as follow:
\bqre 
 R_k^\eps &= &
  \chi_\eps(a_{k1}v_1^0 - b_k,\sum_j a_{kj}
 v_{kj}^1 ) + 
\eps g_k^\eps(s,r_1^\eps,\cdots,r_N^\eps)
 + \eps \alpha_3^\eps (v_k^\eps)'', 
\eqre
where $ g_k^\eps$ is defined by
$
  g_k^\eps(s,r_1^\eps,\cdots,r_N^\eps) 
  = 
 \eps^{-1}\{ ( V_k^\eps  + \eps^2 \sum_j a_{kj}  r_j^\eps )_+
 -
 ( V_k^\eps )_+ \},  
$ 
 and $
V_k^\eps 
  =  a_{k1}v_k^0 -b_k
+ \eps \sum_j a_{kj}
v_{kj}^1. 
$
Notice that $g_k^\eps$ satisfies 
$ \dis
| g_k^\eps(s,r_1^\eps,\cdots,r_N^\eps)  |
\leq  
 \sum_j | a_{kj}| | r_j^\eps |.
$
\\
A key  ingredient is the  energy 
$ 2 E = \dis  \sum_k ( \alpha_0 (r_k')^2 + \lambda_k^2 r_k^2)$
 for the homogeneous system: $ L_k r_k = 0,\; k=1,\cdots,N$ and the 
 for the inhomogeneous system:
 \bqre
  - L_k r_k &=& S_k(s)+ f_k^\eps(s) + \eps g_k(s,r_1^\eps,\cdots,r_N^\eps),
  \eqre
 for $ k=1,\cdots,N$, 
with 
 $
   S_k  = H(a_{k1}v_1^0-b_k) \sum_{j=1}^{N} a_{kj}v^1_j 
            + \alpha_2 v_k^0" + \alpha_1 v_k^1",
$
 and $\alpha_1,\, \alpha_2$ are well chosen to avoid secular term when $k=1$.
 Thus, all $S_k$ are $2\pi$ periodic.
 $S_1$ is not resonant with $L_1$. 
 The $\lambda_k$ are $\zit $ independent. 
 We can apply  Lemma \ref{lemmaTimeSystem}
 which is enough to conclude the proof.
 \epro
\medskip  $\mbox{ }$ \\
%
%
{\bf Proof of Theorem \ref{Th43} :}
 The proof follows two steps. First the solution for linear equations
 \eqref{linsysacc} are bounded by $o(t)$. 
Second, energy estimates are used to bound $r^\eps$.
\\
At the end we prove remark 
                  \ref{polarisation}.
\\
Notice that we do not use  Lemmas \ref{LDLL}, \ref{LHzi}.
Indeed,  we have no term with $\chi_\eps$. 
 We only use that functions  $u_+$ and $v_k^0$ are Lipschitz,
 the Lemma \ref{LSp} to identify resonant terms when $b_k=0$ and an
  energy estimate.  
But, since all modes are excited, the accuracy is weaker 
than the precision obtained in Theorem \ref{Th41}, as in \cite{SP}.
\medskip \\
\underline{ Step 1}: the $N$ problems \eqref{linsysacc} involves 
                decoupled equations rewritten as follow with  $\omega > 0$,
       \bqre 
             \phi''(s)+\omega^2\phi(s) &= & S(s) \in C^0_{ap}(\rit,\rit),  
   \qquad  \pm \omega \nin Sp[S].        
        \eqre 
There is no resonance since $\pm \omega$ are not in the spectrum of $S$.
But, $Sp[S]$ is dense in $\rit$. 
Indeed $\lambda_1,\cdots,\lambda_N$ are $\zit$ independent.
In general, we cannot expect that $\phi$ is bounded on the real line,
 see \cite{Cor}, but $\phi $ is less than $\mathcal{O}(s)$ for large time.
 We can compute explicitly $\phi$
\bqre 
   \phi(s) & = & A \cos(\omega s) + B \sin(\omega s) + \psi(s), \\
   \omega \psi(s) & =& \dis \int_0^s S(\sigma) \sin(\omega(s-\sigma)) d\sigma 
           \\ &  = &
 \dis \sin(\omega s) \int_0^s S(\sigma) \cos(\omega\sigma) d\sigma
                  -\cos(\omega s) \int_0^s S(\sigma) \sin(\omega\sigma) d\sigma
.\eqre 
 The condition $\pm \omega \nin Sp[S]$ is 
 $\dis \lim_{s \rightarrow + \infty} 
 s^{-1}\int_0^s S(\sigma) \exp(\pm i\omega\sigma) d\sigma =0$. 
 That is to say 
\\
 $\dis \int_0^s S(\sigma) \exp(\pm i\omega\sigma) d\sigma =~o(s)$ 
 when $s\rightarrow + \infty$, thus $\psi$ and $\phi$ 
are negligible compared to $s$ for large time. 
\medskip \\
\underline{ Step 2}: Let us decompose the remainder in the following way
 $r_k^\eps = \phi_k + w_k^\eps$. From equation \eqref{RFose} and the previous 
 bound for $\phi_k$ we have in variable $t$ instead of $s_k$ for convenience
 \bqre
      L_k w_k^\eps(t) & = & \mathcal{O}(\eps t) + 
    \left(  \mathcal{O}(\eps \phi_k) + \mathcal{O}(\eps |w^\eps|))
          =  \mathcal{O}(\eps t) + \mathcal{O}(\eps |w^\eps|),
    \right)
 \eqre
 since $\phi_k(t)=o(t).$ Now, we remove the first part of the right hand side 
 with $w_k^\eps= \tilde{w}_k^\eps+z_k^\eps$ and $ \tilde{w}_k^\eps$ is 
 solution of $L_k  \tilde{w}_k^\eps= \mathcal{O}(\eps t)$. 
 Classical energy estimates (or explicit computations as for $\phi$)
 yields to 
   $ \tilde{w}_k^\eps(t)= \mathcal{O}(\eps t^2)$.
 Thus there exists a constant $C_1> 0$ such that $z_k^\eps$ satisfies
 \bqre
  \left | L_k z_k^\eps \right| & \leq &
        C (\eps^2 t^2 + \eps |z^\eps|).
 \eqre
Multiplying each inequality by $|(z_k^\eps)'|$, summing up with respect to $k$, 
integrating on $[0,T]$,  by Cauchy-Schwarz inequality, 
 with $D= 2C (\min(\lambda_k)+ \min(\lambda_k)^2)$ we get
\bqre 
  E(T) & = & \dis \sum_{k=1}^{N}
                    \left( \lambda_1^2 ((z_k^\eps)')^2
                          + \lambda_k^2 (z_k^\eps)^2 
                      \right)\\
     & \leq & 2 C \eps^2 T^2 \dis \int_0^T \sum_{k=1}^{N}|(z_k^\eps)'|(t)dt
               + 2C\eps\dis \int_0^T |(z^\eps)' \cdot z^\eps| dt\\
    & \leq & D \eps^2 T^{2.5} \dis\sqrt{ \int_0^T  E(t)dt}
             + D\eps\dis \int_0^T E(t) dt.
\eqre
Let $Y(T)$ be $\dis\int_0^T E(t) dt $, thus $Y(0)=0$ and for all
 $ t \in [O,T]$, 
\bqre 
   E(t)=Y'(t) &\leq& D \eps^2 T^{2.5} \sqrt{Y(t)} + D\eps Y(t).
\eqre
 Since 
$\dis \int^Y_0 \frac{dy}{A\sqrt{y}+y} =
 2 \ln \left(1 + \frac{\sqrt{y}}{A} \right) $
we obtain $\sqrt{Y(T)} \leq \eps T^{2.5} \exp(D\eps T)$ and 
 then 
 \bqre E(T) &\leq& 2D \eps^3 T^5 \exp(D\eps T).
 \eqre
Finally $r_k^\eps = \phi_k+ \tilde{w}_k^\eps+z_k^\eps
 =  o(T) +\mathcal{O}(\eps T^2)+\mathcal{O}(\eps^{1.5} T^{2.5} \exp(D\eps T)) $,
 so for any $T_\eps = o(\eps^{-1}) $ 
we have in $W^{1,\infty}(0,T_\eps)$ for all $T \leq T_\eps$
\bqre
      \eps r_\eps(T)  &= & o(\eps T_\eps) + \mathcal{O}(\eps^2 T_\eps^2)
                    +\mathcal{O}(\eps^{2.5} T_\eps^{2.5}),                 
\eqre
which is enough  to have the convergence in $W^{1,\infty}(0,T_\eps)$. 
Furthermore $r_k^\eps$ satisfies  the second order differential equation
\eqref{eqNdofRemainderH}  which is enough to get the convergence 
 in $W^{2,\infty}$.
\medskip \\
\underline{ About remark  \ref{polarisation}}:
 From Theorem \ref{Th43}, this result its obvious. 
 Let us explain why we cannot go further up to the order $\eps^2$.

Unfortunately $S_k$ is not periodic since $v_j^1$ is quasi-periodic 
 for $j\neq 1$.
 Indeed, the following  initial conditions
$   
  v_k^1(0)=0, \quad   (v_k^1)'(0)=0,   k \neq 1,
$
yields to  a quasi-periodic function, 
sum of two periodic functions with different periods
$2\pi$ and $2\pi \lambda_1/\lambda_k$, thus a globally bounded function 
$
    v_k^1(s) = \dis  \phi_k^1(s)
   - \phi_k^1(0)\cos\left(\frac{\lambda_k}{\lambda_1}s\right). 
$
 So we cannot apply  
Lemma \ref{lemmaTimeSystem}.
 \\
 Let us decompose $S_k=P_k+Q_k$  for $k\neq 2$ where $P_k$ is periodic and 
  $Q_k$ is almost-periodic 
 \bqre 
   Q_k(s) &=& -H(a_{k1}v_1^0(s)-b_k) \sum_{j=1}^{N} a_{kj}
                \phi^1_j(0)\cos\left(\frac{\lambda_j}{\lambda_1}s\right). 
 \eqre 
Let $w_k$ be a solution of $-L_k w_k = Q_k$ then 
$\dis Sp[w_k] \in \bigcup_{j} 
 \left\{\pm \frac{\lambda_j}{\lambda_1} + \zit  \right\}$, 
 so the spectrum of $w_k$ is discrete  and there is resonance 
 in the $N-1$ equations, 
  $-L_k r^\eps_k=S_k+\cdots $, $k\neq 1$ and the expansion does not still 
 valid for time of the order $\eps^{-1}$.  
 \epro

   %
 

{\it Acknowledgments}: we thank Alain L\'eger,  Vincent Pagneux and 
 St\'ephane Roux 
 for their valuable remarks at the 
fifth meeting of the GDR US,  Anglet, 2008.  
We also  thank G\'erard Iooss for fruitful discussions.

 \bibliographystyle{plain}

 \end{document}